\documentclass[12pt,reqno]{amsart}
1
\usepackage{mathrsfs, amssymb}
\usepackage{array}
\usepackage{tabularx,multirow,makecell,longtable}
\usepackage{mathcomp}
\usepackage{upgreek}
\usepackage[matrix,arrow,curve,all, frame]{xy}
\usepackage{amsthm}
\usepackage[margin=2cm]{geometry}
\usepackage{hyperref}
\usepackage[toc,page]{appendix}
 \usepackage{ extpfeil, extarrows } 
 \usepackage[dvipsnames]{xcolor}
\usepackage{graphicx}
\usepackage{tikz}
\usetikzlibrary{decorations.pathreplacing}
\newenvironment{fulltext}{}{}

\newcommand{\dd}{\mathrm{d}}

\newcommand{\CC}{\mathbb{C}}

\newcommand{\PP}{\mathbb{P}}
\newcommand{\QQ}{\mathbb{Q}}

\newcommand{\mmm}{\mathfrak{m}}
\newcommand{\axis}{\mathrm{\mbox{-}axis}}

\newcommand{\I}{{\mathcal{I}}}
\newcommand{\J}{{\mathcal{J}}}

\newcommand{\MMM}{{\mathscr{M}}}

\newcommand{\OOO}{{\mathscr{O}}}

\newcommand{\LLL}{{\mathscr{L}}}

\newcommand{\di}{\operatorname{div}}
\newcommand{\wt}{\operatorname{wt}}

\newcommand{\Sing}{\operatorname{Sing}}

\newcommand{\gr}{\operatorname{gr}}

\newcommand{\red}{{\operatorname{red}}}

\newcommand{\toplus}{\mathbin{\tilde\oplus}}
\newcommand{\totimes}{\mathbin{\tilde\otimes}}

\newcommand{\F}{{\mathrm{F}}}

\newcommand{\Supp}{{\operatorname{Supp}}}

\newcommand{\mumu}{{\boldsymbol{\mu}}}

\renewcommand{\omega}{\upomega}

\newcommand{\ncomment}[1]{}
\newcommand{\typec}[1]{$\mathrm{\left(#1\right)}$}
\newcommand{\type}[1]{$\mathrm{#1}$}
\newcommand{\typet}[2]{$\mathrm{#1}_{#2}$}
\newcommand{\xref}[1]{\textup{\ref{#1}}}

\renewcommand\labelenumi{\rm (\roman{enumi})}
\renewcommand\theenumi{\rm (\roman{enumi})}

\newcommand{\quo}[1]{#1}
\newcommand{\qqed}{}

\makeatletter
\@addtoreset{equation}{section}
\makeatother

\renewcommand{\theequation}{\arabic{section}.\arabic{equation}}
\newenvironment{zproof}[1][]{\begin{proof}[Proof  #1]}{\end{proof}}
\theoremstyle{plain}
\newtheorem{theorem}{Theorem}[section]
\newtheorem{lemma}[theorem]{Lemma}

\newtheorem{proposition}[theorem]{Proposition}

\newtheorem{corollary}[theorem]{Corollary}

\newtheorem{slemma}[theorem]{Lemma}

\theoremstyle{definition}

\newtheorem{sdefinition}[theorem]{Definition}

\newtheorem{notation}[theorem]{Notation}

\newtheorem{example}[theorem]{Example}

\newtheorem{sremark}[theorem]{Remark}

\title{General elephants \\ for threefold extremal contractions \\ with one-dimensional fibers: exceptional case}
\author{Shigefumi Mori}
\author{Yuri Prokhorov}
 \address{
 S.~Mori:
Kyoto University Institute for Advanced Study,
Kyoto University, Kyoto, Japan
 \newline\indent
Research Institute for Mathematical Sciences,
Kyoto University, Kyoto, Japan
 \newline\indent
Chubu University Academy of Emerging Sciences,
Chubu University, Aichi, Japan
  }
  \email{mori@kurims.kyoto-u.ac.jp}

  \address{Y.~Prokhorov:
  Steklov Mathematical Institute of Russian Academy of Sciences, Moscow, Russian Federation
  }
  \email{prokhoro@mi-ras.ru}

\thanks{ This work was 
supported by the Research Institute for Mathematical
 Sciences, an International  Joint Usage/Research Center located in Kyoto University.
 }

\date{}

\begin{document}
 \maketitle\tableofcontents

 \begin{fulltext}
\begin{abstract}
Let $(X, C)$ be a germ of a threefold $X$ with terminal singularities along a connected reduced
complete curve $C$ with a contraction $f : (X, C) \to (Z, o)$ such that $C = f^{-1} (o)_{\red}$ and $-K_X$ is $f$-ample. Assume that
each irreducible component of $C$ contains at most one point of index $>2$.
We prove that a general member $D\in |{-}K_X|$ is a normal surface with Du Val singularities.
\end{abstract}
\section{Introduction}

The present paper is a continuation of a series of papers on the classification of extremal contractions with one-dimensional fibers (see the survey
\cite{MP-1p} for an introduction).
Recall that an \textit{extremal curve germ} is the analytic germ $(X,\, C)$
of a threefold $X$ with terminal singularities along a reduced connected complete curve $C$ 
such that there exists a contraction $f: (X,\, C)\to (Z,\, o)$ such that $C=f^{-1}(o)_{\red}$ and 
$-K_X$ is $f$-ample. There are three types of extremal curve germs: flipping, divisorial and $\QQ$-conic bundles, and all of them are important building blocks in the three-dimensional minimal model program. 

The first step of the classification is to establish the existence of a ``good'' member of the anti-canonical linear system. This is M. Reid's so-called ``general elephant conjecture'' \cite{Reid:YPG}. In the case of irreducible central curve $C$ the conjecture has been proved:

\begin{theorem}[\quo{\cite[Th.~2.2]{KM92}, \cite{MP:cb3}}]
\label{thm:ge:irr}
Let $(X,\, C)$ be an extremal curve germ with irreducible central curve $C$. Then a general member $D\in |{-}K_X|$ is a normal surface with Du Val singularities. 
\end{theorem}

Moreover, all the possibilities for general members of $|{-}K_X|$ have been classified. Firstly, extremal curve germs with irreducible central curve are divided into two classes: semistable and exceptional. 
Such a germ $(X,\, C)$ is said to be \textit{semistable} if for the restriction of the corresponding contraction $f: (X,\, C)\to (Z,\, o)$ to a general member $D\in |{-}K_X|$, we have the Stein factorization $f_D: D\to D'\to f(D)$ with surface $D'$ having only Du Val singularities of type \type{A} \cite{KM92}. 
Non-semistable extremal curve germs are called \textit{exceptional}. Semistable extremal curve germs are subdivided into two types: \typec{k1A} and \typec{k2A} while exceptional ones are subdivided 
into the following types: \type{cD/2}, \type{cAx/2}, \type{cE/2}, \type{cD/3}, \typec{IIA}, \typec{II^\vee}, \typec{IE^\vee}, \typec{ID^\vee}, \typec{IC}, \typec{IIB}, \typec{kAD}, and \typec{k3A} 
(see \cite{KM92}, \cite{MP:cb1} and \cite{MP:cb3}). \footnote{
This statement is not true as it stands. In the case \typec{IE^\vee}\ $D$ is always of type \type{A_7} and in the case \typec{ID^\vee}\
$D$ is type \type{A_3} or $\mathrm{D}_k$, $k \ge 4$ depending on whether the 
non-Gorenstein point $P \in X$ is of type \type{cA/2} or \type{cAx/2}   \cite[1.2.4]{MP:cb1} (see Appendix~\ref{appB} for the full list of $\QQ$-conic bundle germs with irreducible central 
fiber).

To correct the original statement, the cases \typec{ID^\vee} with $P$ of type \type{cA/2} and  \typec{IE^\vee} should be moved to the group of semistable extremal curve germs, and only the case 
\typec{ID^\vee} with $P$ of type \type{cAx/2} should remain in the group of exceptional ones.}

The result stated in Theorem~\ref{thm:ge:irr} is very important in three-dimensional geometry. For example, the existence of a good member $D\in |{-}K_X|$ for flipping contractions is a sufficient condition for the existence of flips \cite{Kaw:Crep} and the existence of a good member $D\in |{-}K_X|$ in the $\QQ$-conic bundle case proves Iskovskikh's conjecture about singularities of the base \cite{P97}, \cite{MP:cb1}. 

Reid's conjecture also has been proved for arbitrary central curve $C$ in the case of $\QQ$-conic bundles over singular base:

\begin{theorem}[\quo{\cite{MP:cb2}}]
\label{theorem:cb:singular:S}
Let $(X,\, C)$ be a $\QQ$-conic bundle germ and let $f: (X,\, C)\to (Z,\, o)$ be the corresponding contraction. Assume 
that $(Z, o)$ is singular.
Then a general member $D\in |{-}K_X|$ is a normal surface with Du Val singularities.
\end{theorem}

In this paper we study Reid's conjecture for extremal curve germs with reducible central curve. 
Our main result is the following theorem.

\begin{theorem}\label{thm:main}
Let $(X,\, C)$ be an extremal curve germ.
Assume that $(X,\, C)$ satisfies the following condition:
\begin{enumerate}
\renewcommand\labelenumi{\rm (\fnsymbol{enumi})}
\renewcommand\theenumi{\rm (\fnsymbol{enumi})}
\item 
\label{thm:main:condition}
Each irreducible component of $C$ contains at most one point of index $>2$.
\end{enumerate}
Then a general member $D\in |{-}K_X|$ is a normal surface with Du Val singularities.
Moreover, for each irreducible component $C_i\subset C$ with two non-Gorenstein points or of types \typec{IC} or \typec{IIB}, the dual graph $\Delta(D,\, C_i)$
has the same form as the irreducible extremal curve germ $(X,\, C_i)$ \textup(see Theorem~\xref{theorem:ge}\textup).
\end{theorem}

Throughout this paper we use the standard notation \typec{IC}, \typec{IIB} etc for types of extremal curve germs $(X,\, C)$ with irreducible central fiber \cite{KM92}.
Sometimes, to specify the indices of singular points, we will use subscripts. E.g. \typec{kAD_{2,m}} means that the indices of points of $(X,\, C)$ are 2 and $m$.
Some of the subscripts can be omitted if it is not important in the consideration, e.g. \typec{k2A_2} means that $(X,\, C)$ contains a point of index $2$ (and another point of index $>1$). 

According to the classification of birational extremal curve germs 
the condition~\ref{thm:main}\ref{thm:main:condition} is equivalent saying that an arbitrary component $C_i\subset C$ of type \typec{k2A}
has a point of index $2$. 

\begin{corollary}
\label{cor:main}
Let $(X,\, C)$ be an extremal curve germ 
and let $C_i\subset C$ be an irreducible component.
\begin{enumerate}
\item \label{cor:main1}
If $C_i$ is of type \typec{IIB}, then any other component $C_j\subset C$ is of type \typec{IIA} or \typec{II^\vee}.
\item \label{cor:main2}

If $C_i$ is of type \typec{IC} or \typec{k3A}, then any other component $C_j\subset C$ meeting $C_i$ is of type \typec{k1A} or \typec{k2A}.

\item \label{cor:main3}

If $C_i$ is of type \typec{kAD}, then any other component $C_j\subset C$  meeting $C_i$  is of type \typec{k1A}, \typec{k2A}, \type{cD/2}, or \type{cAx/2}.
\item \label{cor:main4}

If $C_i$ is of type \typec{k2A_2}, then any other component $C_j\subset C$  meeting $C_i$  is of type  \typec{k1A}, \typec{IC}, or \typec{k2A_{n,m}}, $n,\, m\ge 3$.
\end{enumerate}
\end{corollary}
There are more restrictions on the combinatorics of the components of $C$.
This will be treated in a subsequent paper.
\ncomment{Examples of extremal curve germs satisfying conditions of Theorem~\ref{thm:main} can be found in the appendix of the arxiv version of this paper.}

\textbf{Acknowledgements.} The paper was written during the second author's visits to
RIMS, Kyoto University. The authors are very grateful to the institute for their
support and hospitality.

\section{Preliminaries}
\subsection{}
Recall that a \emph{contraction} is a proper surjective morphism $f:X\to Z$ of 
normal varieties such that $f_*\OOO_X=\OOO_Z$. 
\begin{sdefinition}
Let $(X,\, C)$ be the analytic germ 
of a threefold with terminal singularities along a reduced connected complete curve. We 
say that $(X,\, C)$ is an \emph{extremal curve germ} if there is a contraction 
$f: (X,\, C)\longrightarrow (Z,\, o)$
such that $C=f^{-1}(o)_{\red}$ and $-K_X$ is $f$-ample. 
Furthermore, $f$ is called \emph{flipping} if its exceptional locus coincides 
with $C$ and \emph{divisorial} if its exceptional locus is two-dimensional. If 
$f$ is not birational, then $Z$ is a surface and $(X,\, C)$ is said to be a 
\emph{$\QQ$-conic bundle germ}.
\end{sdefinition}

\begin{lemma}\label{lemma:components}
Let $(X,\, C)$ be an extremal curve germ. Assume that $C$ is reducible.
Then for any proper connected subcurve $C'\subsetneqq C$ the germ $(X,\, C')$ is a birational extremal curve germ.
\end{lemma}

\begin{proof}
Clearly, there exists a contraction $f': X\to Z'$ of $C'$ over $Z$ \cite[Corollary~1.5]{Mori:flip}. We need to show only that $f'$ is birational.
Assume that $(X,\, C')$ is a $\QQ$-conic bundle germ. Then there exists the following commutative diagram
\begin{equation*}
\xymatrix@R=12pt@C=35pt{
X\ar[d]_f\ar[dr]^{f'} \\
Z& Z'\ar[l]_{\varphi}
}
\end{equation*}
where $f$ and $f'$ are $\QQ$-conic bundles contracting $C$ and $C'$, respectively.
The image $\Gamma:=f'(C'')$ of the remaining part $C'':=C-C'$ is a curve on $Z'$ such that $\varphi(\Gamma)=f(C)$ is a point, say $o\in Z$. 
Hence the fiber $f^{\prime -1}(\Gamma)=f^{-1}(o)$ is two-dimensional, a contradiction. \qqed
\end{proof}

\subsection{}
Recall basic definitions of the $\ell$-structure techniques, see \cite[\S~8]{Mori:flip} for details. 
Let $(X,\, P)$ be three-dimensional terminal singularity of index $m$.
Throughout this paper $\pi: (X^\sharp , P^\sharp ) \to (X, P )$ denotes its index-one cover. For
any object $V$ on $X$ we denote by $V^\sharp$ the pull-back of $V$ on $X^\sharp$.

Let $\LLL$ be a coherent sheaf on $X$ without submodules of
finite length $>0$. An \textit{$\ell$-structure} of $\LLL$ at $P$ is
a coherent sheaf $\LLL^\sharp$ on $X^\sharp$ without submodules of
finite length $>0$ with $\mumu_m$-action endowed with an isomorphism
$(\LLL^\sharp)^{\mumu_m}\simeq \LLL$. An \textit{$\ell$-basis of
$\LLL$ at $P$} is a collection of $\mumu_m$-semi-invariants
$s_1^\sharp,\dots,s_r^\sharp\in \LLL^\sharp$ generating
$\LLL^\sharp$ as an $\OOO_{X^\sharp}$-module at $P^\sharp$. Let $Y$
be a closed subvariety of $X$. Note that
$\LLL$ is an $\OOO_Y$-module if and only if $\LLL^\sharp$ is an
$\OOO_{Y^\sharp}$-module. We say that $\LLL$ is \textit{$\ell$-free}
$\OOO_Y$-module at $P$ if $\LLL^\sharp$ is a free
$\OOO_{Y^\sharp}$-module at $P^\sharp$. If $\LLL$ is $\ell$-free
$\OOO_Y$-module at $P$, then an $\ell$-basis of $\LLL$ at $P$ is
said to be \textit{$\ell$-free} if it is a free
$\OOO_{Y^\sharp}$-basis.

Let $\LLL$ and $\MMM$ be $\OOO_Y$-modules at $P$ with
$\ell$-structures $\LLL\subset \LLL^\sharp$ and $\MMM\subset
\MMM^\sharp$. Define the following operations:
\begin{itemize}
\item
$\LLL\toplus \MMM\subset (\LLL\oplus \MMM)^\sharp$ is an
$\OOO_Y$-module at $P$ with $\ell$-structure
\begin{equation*}
(\LLL\toplus \MMM)^\sharp =\LLL^\sharp\oplus\MMM^\sharp.
\end{equation*}
\item
$\LLL\totimes \MMM\subset (\LLL\otimes \MMM)^\sharp$ is an
$\OOO_Y$-module at $P$ with $\ell$-structure
\begin{equation*}
(\LLL\totimes \MMM)^\sharp
=(\LLL^\sharp\otimes_{\OOO_{X^\sharp}}\MMM^\sharp) /
\operatorname{Sat}_{\LLL^\sharp\otimes\MMM^\sharp}(0),
\end{equation*}
where $\operatorname{Sat}_{\mathcal F_1}\mathcal F_2$ is the saturation of
$\mathcal F_2$ in $\mathcal F_1$.
\end{itemize}
These operations satisfy standard properties (see
\cite[8.8.4]{Mori:flip}). If $X$ is an analytic threefold with
terminal singularities and $Y$ is a closed subscheme of $X$, then
the above local definitions of $\toplus$ and $\totimes$ patch with
corresponding operations on $X\setminus \Sing X$. Therefore, they
give well-defined operations of global $\OOO_Y$-modules.

\begin{lemma}\label{lemma:surfaces}
Let $(D,\, C)$ be the germ of a normal Gorenstein surface along a proper reduced connected curve $C=\cup C_i$, where $C_i$ are irreducible components. 
Assume that the following conditions hold:
\begin{enumerate}
\item\label{lemma:surfaces1}
$K_D\sim 0$,
\item\label{lemma:surfaces2}
there is a birational contraction \mbox{$\varphi: (D, C)\to (R,\, o)$} such that $\varphi^{-1}(o)_\red= C$, 
\item\label{lemma:surfaces3}
there is a point $P\in D$ which is not Du Val of type \type{A}. 
\end{enumerate}
Then $D$ has only Du Val singularities on $C\setminus \{P\}$.
\end{lemma}

\begin{proof}
Assume that there is a point $Q\in D\setminus \{P\}$ which is not Du Val.
If there exists a component $C_i\subset C$ passing through $Q$ but not passing through $P$, we can contract it: $D\to D'$
over $R$.
The contraction is crepant, so the image of $Q$ is again a non-Du Val point. 
Replace $D$ with $D'$.
Continuing the process we may assume that $P$ and $Q$ are connected by some component $C_i\subset C$.
Moreover, by shrinking $C$ we may assume that $C_i=C$, i.e. $C$ is irreducible.
Since $D$ is Gorenstein, the point $Q\in D$ is not log terminal and the point $P\in D$ is log terminal only if it is Du Val of type \type{D} or \type{E}.
Hence the pair $(D, C)$ is not log canonical at $Q$ 
and not purely log terminal at $P$ \cite[Theorem~4.15]{KM:book}. Let $H$ be a general hyperplane section passing through $P$.
For some $0<\epsilon,\, \delta \ll1$ the pair $\bigl(D, (1-\epsilon )C+\delta H\bigr)$ is not log canonical at $P$ and $Q$.
Since $-\bigl(K_D+(1-\epsilon )C+\delta H\bigr)$ is $\varphi$-ample, this contradicts Shokurov's connectedness lemma 
\cite{Shokurov:3flip}. \qqed
\end{proof}

\section{Low index cases}
Extremal curve germs of index 2 with arbitrary central curve have been completely classified in \cite[\S~4]{KM92} and 
\cite[\S~12]{MP:cb1}. As an easy consequence we have the following.

\begin{proposition}\label{proposition:index2}
Let $(X,\, C)$ be an extremal curve germ. Assume that all the singularities of $X$ are of index $1$ or $2$, that is, $2K_X$ is Cartier. Then 
a general member $D\in |{-}K_X|$
is a normal surface with Du Val singularities and $D$ does not contain any component of $C$.
\end{proposition}

\begin{proof}
Since the case where $X$ is Gorenstein is trivial, we assume that $X$ has at least one point, say $P$, of index $2$.
In the birational case there are no other non-Gorenstein points and all the components $C_i\subset C$ pass through $P$ \cite[Prop.~4.6]{KM92}.
By \cite[Th.~2.2]{KM92} a general local member $D\in |{-}K_{(X,\, P)}|$ is in fact a general member of $|{-}K_X|$ and this $D$ has only Du Val singularity (at $P$) \cite[(6.3)]{Reid:YPG}.
For the $\QQ$-conic bundle case we refer to \cite[Proof of~12.1]{MP:cb1} and \cite[Corollary~1.4]{MP:cb2}.\qqed
\end{proof}

\begin{proposition}\label{proposition:super-easy-cases}
Let $(X,\, C)$ be an extremal curve germ. Assume that $C$ is reducible and $(X,\, C)$ contains a point $P$ of one of the types \type{cD/2}, \type{cAx/2}, \type{cE/2}, \type{cD/3}.
Then one of the following holds.
\begin{enumerate}
\item \label{proposition:super-easy-cases1}
$P$ is the only non-Gorenstein point of $X$, all the components pass through $P$ 
and do not meet each other elsewhere, and a general member $D\in |{-}K_X|$
is a normal surface with Du Val singularities. Moreover, $D\cap C=\{P\}$.
\item \label{proposition:super-easy-cases2}
There is a component $C_i\subset C$ passing through $P$ such that the germ $(X,\, C_i)$ is divisorial of type \typec{kAD}. Moreover, $(X,\, P)$ is a singularity of type \type{cD/2} or \type{cAx/2}.
\end{enumerate}
\end{proposition}

\begin{proof}
Recall that the intersection points $C_i\cap C_j$ of different components $C_i, C_j\subset C$ 
are non-Gorenstein by \cite[Corollary~1.15]{Mori:flip}, \cite[Prop.~4.2]{Kollar-1999-R} 
and also by \cite[Lemma~4.4.2]{MP:cb1}. 
If $P$ is the only non-Gorenstein point of $X$, then a general member $D\in |{-}K_{(X,\, P)}|$ is in fact a general member of $|{-}K_X|$ \cite[(0.4.14)]{Mori:flip}. This $D$ has only Du Val singularity (at $P$) \cite[(6.3)]{Reid:YPG}. If there exists a non-Gorenstein point $Q\in X$ other than $P$, then we may assume that $Q$ lies on some component $C_i\subset C$ passing through $P$. Thus $(X,\, C_i)$ is a birational extremal curve germ with two non-Gorenstein points (see Lemma~\ref{lemma:components}).
According to \cite[Th.~2.2]{KM92} and \cite{Mori-err} the germ $(X,\, C_i)$ is divisorial of type \typec{kAD} and $(X,\, P)$ is a singularity of type \type{cD/2} or \type{cAx/2}.
This proves the proposition.\qqed
\end{proof}

\section{Extension techniques}

\begin{theorem}[\quo{\cite[Th.~7.3]{Mori:flip}, \cite[Prop.~1.3.7]{MP:cb1}}]
Let $(X,\, C\simeq \PP^1)$ be an irreducible extremal curve germ satisfying the condition~\xref{thm:main}\xref{thm:main:condition}.
Then for a general member $S\in |-2K_X|$ one has $S\cap C=\{P\}$, where $P$ is the point of index $r>2$ or a smooth point \textup(if $(X,\, C)$ is of index $2$\textup). 
Moreover, the pair $\left(X,\, \frac12 S\right)$ is log terminal.
\end{theorem}

\begin{proposition}[\quo{{\cite[Lemma~2.5]{KM92},} {\cite[Prop.~2.1]{MP:cb3}}}]
\label{prop:extension}
Let $(X, C)$ be an extremal curve germ \textup($C$ is not necessarily irreducible\textup) and let $S\in |-2K_X|$ be a general member.
Assume that the set $\Sigma:=S\cap C$ is finite.

\begin{enumerate}
\item \label{prop:extension:bir}
If $(X, C)$ is birational, then the natural map 
\begin{equation}\label{prop:extension:cb:delta}
\tau: H^0\bigl(X,\, \OOO_X(-K_X)\bigr)\xlongrightarrow{\hspace*{17pt}} \omega_{(S,\, \Sigma)}= H^0\bigl(S,\, \OOO_S(-K_X)\bigr)
\end{equation}
is surjective, where $\omega_{(S,\, \Sigma)}$ is the dualizing sheaf of $(S,\, \Sigma)$.

\item \label{prop:extension:cb:delta-bar}
If $(X, C)$ is a $\QQ$-conic bundle germ over a smooth base surface, then the natural map 
\begin{equation}\label{prop:extension:cb:eq}
\bar \tau: H^0\bigl(X,\, \OOO_X(-K_X)\bigr)\xlongrightarrow{\hspace*{17pt}} \omega_{(S,\, \Sigma)}/\Omega^2_{(S,\, \Sigma)}
\end{equation}
is surjective, where $\Omega^2_{(S,\, \Sigma)}$ is the sheaf of holomorphic $2$-forms on $(S,\, \Sigma)$.
\item \label{prop:extension:cb-a}
If $(X, C)$ is a $\QQ$-conic bundle germ over a smooth base surface and $\Sigma=\Sigma_1 \amalg\Sigma_2$, $\Sigma_i\neq\varnothing$, then
\begin{equation}\label{prop:extension:delta1}
\tau_1: H^0\bigl(X,\, \OOO_X(-K_X)\bigr)\xlongrightarrow{\hspace*{17pt}} \omega_{(S,\, \Sigma_1)}
\end{equation}
is surjective.
\end{enumerate}
\end{proposition}

\begin{proof}
For~\ref {prop:extension:bir} we refer to \cite[Lemma~2.5]{KM92}.
Let us show~\ref {prop:extension:cb:delta-bar}.
Note that by the adjunction $\OOO_S(K_S)= \OOO_S(-K_X)$. Let $f:(X,\, C)\to (Z,\, o)$ be the corresponding $\QQ$-conic bundle contraction and let $g=f|_S: S\to Z$ be its restriction to $S$. Since
the base surface $Z$ is smooth, by \cite[Lemma~4.1]{MP:cb1} there is a canonical isomorphism
\begin{equation*}
R^1f^* \omega_X \simeq \omega_Z.
\end{equation*}
Then we apply \cite[Prop.~2.1]{MP:cb3} in our situation:
\begin{equation*}
\xymatrix@R=9pt{
H^0\bigl(X,\, \OOO_X(-K_X)\bigr)\ar[r]\ar@{=}[d] &H^0\bigl(S,\, \OOO_S(-K_X)\bigr)\ar[d]^{\rotatebox{-90}{$\simeq$}} 
\\
f_*\omega_X(S)\ar[r]\ar@{->>}[rd] & \omega_{(S,\, \Sigma)}\ar[d]
\\
& \omega_{(S,\, \Sigma)}/ g^*\omega_{(Z,\, o)}\ar@{->>}[r] &\omega_{(S,\, \Sigma)}/\Omega^2_{(S,\, \Sigma)}
}
\end{equation*}
and obtain the surjectivity of $\tau$.

For~\ref{prop:extension:cb-a} we consider the map $g_i: S_i\to Z$ which is the restriction of $g$ to $S_i=(S,\, \Sigma_i)\subset S$ and the induced exact sequence
\begin{equation*}
\xymatrix@R=7pt{
f_*\omega_X(S) \ar[r] & \omega_{(S,\, \Sigma_1)}\oplus \omega_{(S,\, \Sigma_2)}\ar[r]& \omega_Z\ar[r]\ar@/^15pt/[l]^{g_2^*} & 0 
}
\end{equation*}
Then we see that $g_2^*: \omega_Z\to 0\oplus \omega_{(S,\, \Sigma_2)}$ is a splitting homomorphism. Therefore, the homomorphism
\begin{equation*}
f_*\omega_X(S)\xlongrightarrow{\hspace*{17pt}} \omega_{(S,\, \Sigma_1)}\oplus\left(\omega_{(S,\, \Sigma_2)}/ g_2^*\omega_{(Z,0)}\right)
\end{equation*}
is surjective. \qqed
\end{proof}

\begin{lemma}\label{lemma:index2point}
Let $(\bar{X},\, \bar{C})$ be an extremal curve germ with reducible central curve $\bar{C}$.
Suppose $\bar X$ satisfy the condition~\xref{thm:main}\xref{thm:main:condition} and that there is a component $C\subset \bar C$ of type \typec{k1A} which meets $\bar C-C$ at a point $P$ of index $2$. Then a general member $D\in |{-}K_{\bar X}|$ does not contain $C$.
\end{lemma}

\begin{proof}
On each irreducible component $C_i$ of $\bar{C}$ there exist at most one point of index $>2$. Let $\{P_a\}_{a \in A}$ be the collection of such points.
For each $C_i$ without points of index $>2$, choose
one general point of $C_i$. Let $\{P_b\}_{b \in B}$ be
the collection of such points.
For each $i \in A \cup B$, let $S_i \in |-2K_{(X,\, P_i)}|$
be a general element on the germ $(X,\, P_i)$, and set $S=\sum_{i \in A \cup B} S_i$. Then
$S$ extends to an element $|-2K_X|$ by \cite[Thm. (7.3)]{Mori:flip}.
A generator $\sigma_b$ of $\OOO_{S,b}(-K_X)\simeq \OOO_{S,b}$ lifts to $s \in H^0\bigl(X,\, \OOO_X(-K_X)\bigr)$ by~\ref{prop:extension}.\ref{prop:extension:bir} if $(\bar X, \bar C)$ is birational and by $A\neq\varnothing$ and~\ref{prop:extension}.\ref{prop:extension:cb:delta-bar} if otherwise. In either case we have $C \not\subset D$. \qqed
\end{proof}

\section{Review of \cite[\S~2]{KM92}} 
We need some refinement of facts on birational extremal curve germs with irreducible central fiber proved in \cite[\S~2]{KM92}.

\subsection{}
\label{not:graphs}
Below, for a normal surface $D$ and a curve $C\subset D$,
we use the usual notation of graphs $\Delta (D,\, C)$
of the minimal resolution of $D$ near $C$:
each vertex labeled $\bullet$ corresponds to an irreducible component of $C$ and each
$\circ$ corresponds to a component $E_i\subset E$ of the  exceptional divisor $E$ on the minimal
resolution of $D$. 
Note that in our situation below  $E_i^2=-2$ for all $E_i$.

\begin{theorem}[\quo{\cite[Th.~2.2]{KM92}, \cite{Mori-err}}]
\label{theorem:ge}
Let $(X,\, C\simeq\PP^1)$ be a birational extremal curve germ and let $D\in |{-}K_X|$ be a general member. Then $D$ is a normal surface with Du Val singularities.
Moreover, either $D\cap C$ is a point or $D\supset C$ and for the graph $\Delta(D,\, C)$ one and only one of the following possibilities holds:
\noindent
\renewcommand{\arraystretch}{1.1}
\setlength{\tabcolsep}{12pt}
\begin{longtable}{p{47pt}l}
\typec{IC} 
\newline $\scriptstyle{k=1}$
& $\xymatrix@R=3pt@C=17pt{
{\underbrace{\circ -\cdots - \circ}_{m-3\ge 2}} \ar@{-}[r]&\circ\ar@{-}[r]\ar@{-}[d]&\circ
\\
&\bullet
}$
\\
\typec{IIB} \newline $\scriptstyle k=3$, $\scriptstyle m=4$&
$\xymatrix@R=6pt@C=17pt{
\circ\ar@{-}[r]&\circ\ar@{-}[r]&\circ\ar@{-}[r]&\circ\ar@{-}[r]&\bullet 
\\
&&\circ\ar@{-}[u]
}$  
\\
\typec{kAD} \newline  $\scriptstyle k=1$,\ $\scriptstyle n=2$ 
&
\raisebox{11pt}{$\xymatrix@R=-4pt@C=17pt{
&&&&\circ
\\
{\underbrace{\circ -\cdots - \circ}_{m-1\ge 2}} \ar@{-}[r]&{\bullet}\ar@{-}[r]&{\underbrace{\circ -\cdots - \circ}_{2l-2\ge 0}} \ar@{-}[r]&\circ \ar@{-}[ru]\ar@{-}[rd]
\\
&&&&\circ
}$ }
\\
\typec{k3A} \newline  $\scriptstyle k=1$,\ $\scriptstyle n=2$  
& 
\raisebox{11pt}{$\xymatrix@R=-4pt@C=17pt{
&&\circ
\\
{\underbrace{\circ -\cdots - \circ}_{m-1\ge 2}} \ar@{-}[r]&\bullet \ar@{-}[ru]\ar@{-}[rd]
\\
&&\circ
} $ }
\\
\typec{k2A} &
$
\xymatrix@R=3pt@C=17pt{
{\underbrace{\circ -\cdots - \circ}_{km-1}} \ar@{-}[r]&\bullet\ar@{-}[r]&{\underbrace{\circ -\cdots - \circ}_{ln-1}}\ar@{}[r]&
}$
\end{longtable}
\noindent
where $m$ and $k$ are the index and axial multiplicity \cite[1a.5(iii)]{Mori:flip} of a singular point of $X$, and $n$ and $l$ are the ones for the other non-Gorenstein point
\textup(if any\textup).
\end{theorem}

In the cases \typec{IC}, \typec{IIB}, \typec{kAD}, \typec{k3A}, and \typec{k2A_2} Theorem~\ref{theorem:ge} is a consequence of the following.

\begin{theorem}[\quo{cf. {\cite[\S~2]{KM92}}, \cite{Mori-err}}]
\label{theorem:extending}
Let $(X,\, C)$ be a birational extremal curve germ with irreducible central curve of type \typec{IC}, \typec{IIB}, \typec{kAD}, \typec{k3A}, or \typec{k2A_2}.
Let $S\in |-2K_X|$ be a general member \textup(so that $S\cap C=\{P\}$, where $P$ is the point of index $r>2$\textup). Let $\sigma_S\in H^0\bigl(S, \OOO_S(-K_X)\bigr)$ be a general section.
Then for any section $\sigma\in H^0\bigl(X,\, \OOO_X(-K_X)\bigr)$ such that 
\begin{equation}
\label{eq:sigma-sigma}
\sigma|_S\equiv \sigma_S\mod \Omega^2_S, 
\end{equation}
\textup(see~\eqref{prop:extension:cb:eq}\textup) the divisor $D:=\di(\sigma)$ is a normal surface with only Du Val singularities.
Furthermore, the configuration $\Delta(D,\, C)$ is as described in Theorem~\xref{theorem:ge}. 
\end{theorem}

Below we outline the proof of Theorem~\ref{theorem:extending} following \cite[\S~2]{KM92}.
We treat the possibilities \typec{IC}, \typec{IIB}, \typec{k3A}, \typec{kAD}, and \typec{k2A_2} case by case.

\subsection{Case \typec{IC}}
By \cite[(A.3)]{Mori:flip} we have the following identification at $P$
\begin{equation*}
(X,\, C)= \left(\CC^3_{y_1,\, y_2,\, y_4},\, \ \{y_1^{m-2}-y_2^2= y_4=0 \} \right) /\mumu_{m}(2,m-2,1).
\end{equation*}
A general divisor $S\in |-2K_X|$ is given by $y_1=\xi(y_2,\, y_4)$ with $\xi \in (y_2,\, y_4)^2$ such that $\wt(\xi) \equiv 2 \mod m$. Thus we have
\begin{gather}
S\simeq \CC^2_{y_2,\, y_4} /\mumu_{m}(m-2,1),\quad 
\omega_S=(\OOO_{S^\sharp,\, P^\sharp}\, \dd y_2\wedge \dd y_4)^{\mumu_m},
\\
\label{eq:IC:omegaCC}
\omega_S \otimes \CC_P = \CC \cdot y_2^{(m-1)/2}\, \dd y_2\wedge \dd y_4\oplus \CC\cdot y_4\, \dd y_2\wedge \dd y_4. 
\end{gather}
Furthermore, 
\begin{equation}
\label{eq:new-new:gr0}
\gr_C^0\omega^*=(P^\sharp)=\left(-1+\textstyle{\frac{m+1}2}\cdot 2P^\sharp\right)\simeq \OOO_C(-1),
\end{equation}
where
\begin{equation*}
\Omega^{-1}:=(\dd y_1\wedge \dd y_2\wedge \dd y_4)^{-1}
\end{equation*}
is an $\ell$-free $\ell$-basis at $P$. Hence, $H^0(C,\, \gr_C^0 \omega^*)=0$ and 
\begin{equation}
\label{eq:IC:H0}
H^0\bigl(X,\, \OOO_X(-K_X)\bigr)=H^0\bigl(X,\, \I_C \totimes\OOO_X(-K_X)\bigr),
\end{equation}
where $\I_C$ is the defining ideal of $C$ in $X$.
Furthermore, by \cite[(2.10.4)]{KM92}
\begin{equation}
\label{eq:new-new:gr1}
\gr_C^1\omega^*=(5P^\sharp)\toplus (0),
\end{equation}
where the $\mumu_m$-semi-invariants
\begin{equation}
\label{eq:IC:basis}
(y_1^{m-2}-y_2^2)\cdot \Omega^{-1},\qquad y_4\cdot \Omega^{-1}
\end{equation}
form an $\ell$-free $\ell$-basis at $P$. Therefore,
\begin{equation*}
\gr_C^1\omega^*\simeq
\begin{cases}
\OOO_C(-1)\oplus \OOO_C & \text{if $m\ge 9$,}
\\
\OOO_C\oplus \OOO_C & \text{if $m=7$,}
\\
\OOO_C(1)\oplus \OOO_C & \text{if $m=5$.}
\end{cases}
\end{equation*}
we have natural homomorphisms 
\begin{equation*}
\delta : 
H^0\bigl(X,\, \OOO_X(-K_X)\bigr) \xlongrightarrow{\hspace*{17pt}} \gr_C^1 \omega^*
\xlongrightarrow{\hspace*{17pt}} (\gr_C^1 \omega^*)^\sharp \otimes \CC_{P^\sharp}.
\end{equation*}
Since $(y_1-\xi)\cdot (\gr_C^1 \omega^*)^\sharp \otimes \CC_{P^\sharp}=0$, the map $\delta$ factors as
\begin{equation*}
\delta :
H^0\bigl(X,\, \OOO_X(-K_X)\bigr) \xlongrightarrow{\hspace*{17pt}} \omega_S \xlongrightarrow{\hspace*{17pt}} (\gr_C^1 \omega^*)^\sharp \otimes \CC_{P^\sharp}.
\end{equation*}
As in \cite[(3.1.1)]{MP:cb3} we see
\begin{equation*}
\Omega^2_S 
\subset \bigl(\mathfrak m_{S,\, P}\cdot y_4 +\mathfrak m_{S,\, P}\cdot y_2^{(m-1)/2}\bigr) \dd y_2\wedge \dd y_4 =\mathfrak m_{S,\, P}\cdot \omega_S,
\end{equation*}
because for arbitrary elements $\phi_1$, $\phi_2$ of the set of generators 
\begin{equation*}
\left\{y_2^m,\, \ y_4^m,\, \ y_2y_4^2,\, \ y_2^{(m-1)/2}y_4\right\}
\end{equation*}
of the ring $\OOO_{S^\sharp}^{\mumu_m}$ we have
\begin{equation*}
\dd\phi_1 \wedge \dd\phi_2 \in \bigl((y_2,\, y_4)y_4+(y_2,\, y_4)y_2^{(m-1)/2}\bigr)
\dd y_2\wedge \dd y_4.
\end{equation*}
Thus $\delta$ factors further as follows:
\begin{equation*}
\delta :
H^0\bigl(\OOO_X(-K_X)\bigr) 
\xtwoheadrightarrow{\hspace*{10pt}} 
\omega_{S}/\Omega^2_{S} 
\xtwoheadrightarrow{\hspace*{10pt}} 
\omega_S \otimes \CC_P 
\xlongrightarrow{\hspace*{17pt}} 
(\gr_C^1 \omega^*)^\sharp \otimes \CC_{P^\sharp},
\end{equation*}
where the last map is a surjection if $m=5$ and the image is generated by $y_4\Omega^{-1}$ if $m\ge 7$ (see \eqref{eq:IC:omegaCC} and \eqref{eq:IC:basis}).
If $m\ge 7$, this implies that the coefficient of $y_4\Omega^{-1}$ in $\sigma_S$ is nonzero. 
If $m=5$, then the coefficients of $y_4\Omega^{-1}$ and $(y_1^{m-2}-y_2^2)\Omega^{-1}$ in $\sigma_S$ are independent and the image $\bar \sigma$ of $\sigma$ in $\gr_C^1\omega^*$ is not contained in $\OOO_C(1)$. Hence, $\bar \sigma$ is nowhere vanishing and so the singular locus of $D$ does not meet $C\setminus \{P\}$. Then we again can take $\sigma_S$ so that it contains the term $y_4\Omega^{-1}$. 
Therefore, $D\in |{-}K_X|$ can be given by the equation $y_4+\cdots =0$. 
Then \cite[Computation~2.10.5]{KM92} shows that $D$ is Du Val at $P$ and that its graph is
as given as type \typec{IC} in~\ref{theorem:ge}.

\subsection{Case \typec{IIB}}
Then, by \cite[(A.3)]{Mori:flip}, the germ $(X,\, C)$ at $P$ can be given as follows
\begin{gather*}
(X,\, C) = \left(\{\phi=0\} {\subset} \CC^4_{y_1,\dots,y_4},\, \{y_1^2-y_2^3=y_3=y_4=0\}\right)/\mumu_4(3,2,1,1),
\\
\phi=y_1^2-y_2^3+\psi, \quad 
\wt(\psi)\equiv 2\mod 4,\quad \psi (0,0,\, y_3,\, y_4) \notin 
(y_3,\, y_4)^3. 
\end{gather*}

A general divisor $S\in |-2K_X|$ is given by $y_2=\xi(y_1,\, y_3,\, y_4)$ with $\xi \in (y_1,\, y_3,\, y_4)^2$ such that $\wt(\xi) \equiv 2 \mod 4$. Thus $S$ is the quotient of a hypersurface $\phi(y_1,\, \xi,\, y_3,\, y_4)=0$ in $\CC^3_{y_1,\, y_3,\, y_4}$ by $\mumu_4(3,1,1)$. We have
\begin{eqnarray}
\omega_S&=&\left(\OOO_{S^\sharp,\, P^\sharp}\, \frac{\dd y_3\wedge \dd y_4}{y_1+\cdots}\right)^{\mumu_4}\nonumber 
\\
\label{eq:IIB:omegaCC}
\omega_S \otimes \CC_P &=& \CC \cdot y_3 \frac{\dd y_3\wedge \dd y_4}{y_1 + \cdots}
\oplus \CC\cdot y_4 \frac{\dd y_3\wedge \dd y_4}{y_1 + \cdots}.
\end{eqnarray}
Furthermore, 
\begin{equation}
\label{eq:new-new:gr0-star}
\gr_C^0\omega^*=(P^\sharp)=\left(-1+3P^\sharp+ 2P^\sharp\right)\simeq \OOO_C(-1),
\end{equation}
where
\begin{equation*}
\Omega^{-1}:=\left(\frac{\dd y_2\wedge \dd y_3\wedge \dd y_4}{\partial \phi/\partial y_1}\right)^{-1}
\end{equation*}
is an $\ell$-free $\ell$-basis at $P$. Hence, $H^0(C,\, \gr_C^0 \omega^*)=0$ and 
\begin{equation}
\label{eq:IIB:H0}
H^0\bigl(X,\, \OOO_X(-K_X)\bigr)=H^0\bigl(X,\, \I_C \totimes\OOO_X(-K_X)\bigr)
\end{equation}
where $\I_C$ is the defining ideal of $C$ in $X$.
Furthermore, by \cite[(2.11)]{KM92}
\begin{equation}
\label{eq:new-new:gr1-star}
\gr_C^1\omega^*=(0)\toplus (1)\simeq \OOO_C\oplus\OOO_C(1),
\end{equation}
where the $\mumu_m$-invariants
\begin{equation}
\label{eq:IIB:basis}
y_3\cdot \Omega^{-1},\qquad y_4\cdot \Omega^{-1}
\end{equation}
form an $\ell$-free $\ell$-basis at $P$. 
As in \typec{IC}, we have natural homomorphisms
\begin{equation*}
\delta : H^0\bigl(\OOO_X(-K_X)\bigr) 
\xtwoheadrightarrow{\hspace*{10pt}} 
\omega_S \otimes \CC_P  
\xlongrightarrow{\hspace*{17pt}}
(\gr_C^1 \omega^*)^\sharp \otimes \CC_{P^\sharp},
\end{equation*}
where the last homomorphism is an isomorphism (see \eqref{eq:IIB:omegaCC} and \eqref{eq:IIB:basis}).
Thus the coefficients of $y_3 \Omega^{-1}$ and $y_4 \Omega^{-1}$ in $\sigma_S$ are independent, whence \cite[Computation~2.11.2]{KM92} shows that $D$ is Du Val at $P$, 
and the image $\bar \sigma$ of $\sigma$ in $\gr_C^1 \omega^*$ is not contained in $\OOO_C(1)$. Hence $\bar \sigma$ is nowhere vanishing and $D$ is smooth outside $P$. Hence the graph $\Delta(D,\, C)$ is as
given as type \typec{IIB} in~\ref{theorem:ge}.

\subsection{Case \typec{k3A}}
\label{case:k3A}
The configuration of singular points on
$(X,\, C)$ is the following: 
a type~\typec{IA} point $P$ of odd index $m\ge 3$, a type~\typec{IA} point $Q$ of index $2$ and a type~\typec{III} point $R$.
According to \cite[(A.3)]{Mori:flip} and \cite[(2.12)]{KM92} we can express 
\begin{equation*}
\begin{array}{ll}
(X,\, C,\, P)&=\left(\CC^3_{y_1,\, y_2,\, y_3},\, (y_1\axis),\, 0\right) /\mumu_{m}(1,(m+1)/2,-1),
\\[3pt] 
(X,\, C,\, Q)&=\left(\CC^3_{z_1,\, z_2,\, z_3},\, (z_1\axis),\, 0 \right)/\mumu_2(1,1,1), 
\\[3pt] 
(X,\, C,\, R)&=\bigl(\{\gamma(w_1,\, w_2,\, w_3,\, w_4)=0\},\, (w_1\axis),\, 0\bigr), 
\end{array}
\end{equation*} 
where $\gamma \equiv w_1w_3 \mod (w_2,\, w_3,\, w_4)^2$. 

For a general divisor $S\in |-2K_X|$ we have $S\cap C=\{P\}$ and $S$ is given by $y_1=\xi(y_2,\, y_3)$ with $\xi \in (y_2,\, y_3)^2$ such that $\wt(\xi )\equiv 1 \mod m$. Thus 
\begin{gather}
S\simeq \CC_{y_2,\, y_3}^2/\mumu_m\left({\textstyle \frac {m+1}2},-1\right),\quad 
\omega_S=(\OOO_{S^\sharp,\, P^\sharp}\, \dd y_2\wedge \dd y_3)^{\mumu_m},
\\
\label{eq:k3A:omegaCC}
\omega_S \otimes \CC_P = \CC \cdot y_2\, \dd y_2\wedge \dd y_3\oplus \CC\cdot y_3^{(m-1)/2}\, \dd y_2\wedge \dd y_3. 
\end{gather}
By the proof of \cite[(2.12.2)]{KM92} we have 
\begin{equation}
\label{eq:new-new:gr0-a}
\gr_C^0\omega^*=\left(-1+{\textstyle \frac {m+1}2}P^\sharp+ Q^\sharp\right)\simeq \OOO_C(-1),
\end{equation}
where an $\ell$-free $\ell$-basis at $P$, $Q$ and $R$, respectively, can be written as follows: 
\begin{eqnarray*}
\Omega_P^{-1}&:=&\left(\dd y_1\wedge \dd y_2\wedge \dd y_3\right)^{-1}
\\
\Omega_Q^{-1}&:=&\left(\dd z_1\wedge \dd z_2\wedge \dd z_3\right)^{-1}
\\
\Omega_R^{-1}&:=&\left(\frac{\dd w_2\wedge \dd w_3\wedge \dd w_4}{\partial \gamma/\partial w_1}\right)^{-1}
\end{eqnarray*}
Hence, $H^0(C,\, \gr_C^0 \omega^*)=0$ and 
\begin{equation}
\label{eq:k3A:H0}
H^0\left(X,\, \OOO_X(-K_X)\right)=H^0\left(X,\, \I_C \totimes\OOO_X(-K_X)\right),
\end{equation}
where $\I_C$ is the defining ideal of $C$ in $X$.
Furthermore, as in \cite[(2.12.4)]{KM92}
we can further arrange that
\begin{equation}
\label{eq:new-new:gr1-a}
\gr_C^1\omega^*=(0)\toplus\left (-1+{\textstyle \frac {m+3}2}P^\sharp\right),
\end{equation}
where $y_2\cdot \Omega_P^{-1}$,\, $z_2\cdot \Omega_Q^{-1}$, $w_2\cdot \Omega_R^{-1}$ 
form an $\ell$-free $\ell$-basis for $(0)$ 
at $P$, $Q$, and $R$, respectively, and $y_3\cdot \Omega_P^{-1}$,\, 
$z_3\cdot \Omega_Q^{-1}$,\, $w_4\cdot \Omega_R^{-1}$
for $(-1+(m+3)/2P^\sharp)$. Moreover,
\begin{equation*}
\gamma \equiv w_1w_3+c_1w_4^2+c_2w_4w_2+c_3w_2^2 \mod (w_3,\, w_2^2,\, w_2w_4,\, w_4^2)\cdot \I_C
\end{equation*}
for some $c_1,\, c_2,\, c_3 \in \CC$ such that $c_1 \neq0$ if $m \ge 5$ 
(see \cite[(2.12.6)]{KM92} and \cite[Remark~2]{Mori-err})
and $(c_1,\, c_2,\, c_3)\neq0$ if $m=3$ (see \cite[(2.12.7)]{KM92} and \cite[Remark~2]{Mori-err}).

As in \cite[(3.1.1)]{MP:cb3} we see
\begin{equation*}
\Omega^2_S 
\subset \left(\mathfrak m_{S,\, P}\cdot y_2 +\mathfrak m_{S,\, P}\cdot y_3^{(m-1)/2}\right) \dd y_2\wedge \dd y_3 =\mathfrak m_{S,\, P}\cdot \omega_S,
\end{equation*}
because for arbitrary elements $\phi_1$, $\phi_2$ of the set of generators 
\begin{equation*}
\left\{y_2^m,\, \ y_3^m,\, \ y_2^2y_3,\, \ y_2y_3^{(m+1)/2}\right\}
\end{equation*}
of the ring $\OOO_{S^\sharp}^{\mumu_m}$ we have
\begin{equation*}
\dd\phi_1 \wedge \dd\phi_2 \in \left((y_2,\, y_3)y_2+(y_2,\, y_3)y_3^{(m-1)/2}\right)
\dd y_2\wedge \dd y_3.
\end{equation*}
Thus the image of the homomorphism
\begin{equation*}
\delta : H^0\bigl(\OOO_X(-K_X)\bigr) 
\xtwoheadrightarrow{\hspace*{10pt}} 
\omega_S \otimes \CC_P 
\xlongrightarrow{\hspace*{17pt}}
\left(\gr_C^1 \omega^*\right)^\sharp \otimes \CC_{P^\sharp},
\end{equation*}
is equal to $(\gr_C^1 \omega^*)^\sharp \otimes \CC_{P^\sharp}$ if $m=3$, and $\CC\cdot y_2\Omega_P^{-1}$ if $m \ge 5$.

If $m\ge 5$, this implies that the coefficient of $y_2\Omega_P^{-1}$ in the image $\bar \sigma$ of $\sigma$ in
$\gr_C^1 \omega^*$ is nonzero and hence nowhere vanishing.
If $m=3$, then the coefficients of $y_2\Omega_P^{-1}$ and $y_3\Omega_P^{-1}$ are independent and hence $\bar \sigma$ is a general global section of $\gr_C^1 \omega^* \simeq \OOO_C\oplus\OOO_C$.
Then the proof of \cite[(2.12.5)]{KM92} shows that $D$ is Du Val and that its graph is as given as type \typec{k3A} in~\ref{theorem:ge}.

\begin{slemma}\label{new-lemma:2}
In the notation of \xref{case:k3A} there exists a deformation $(X_\lambda, C_\lambda\simeq \PP^1)$ of $(X,C)$ which is trivial outside $R$ such that  for $\lambda\neq 0$ the germ $(X_\lambda, C_\lambda)$ 
has a cyclic quotient singularity at $Q$, and is of type \typec{kAD} of case \eqref{scase:kAD:1} (resp. \typec{k2A_2}) if $m \ge 5$ (resp. $m=3$).
\end{slemma}

\begin{proof}
Let $(X_\lambda, C_\lambda)$  be  the twisted extension \cite[1b.8.1]{Mori:flip} of the germ
\begin{equation*}
(X_\lambda, R) =\{ \gamma- \lambda w_2=0\}\supset (C_\lambda, R)= (w_1\axis)
\end{equation*}
by $u=(w_2, w_4)$. Then in $\gr_{C_\lambda}^1 \OOO$ we have  $w_1w_3=\lambda w_2$
 for $\lambda\neq 0$. Since $\gr^1_C\upomega^*=\OOO_C \cdot w_2\Omega^{-1}_R\oplus \OOO_C\cdot w_4 \Omega^{-1}_R$
at $R$, we have at $R$:
\begin{equation*}
 \gr^1_{C_\lambda}\upomega^*=
  \OOO_{C_\lambda} \cdot w_3\Omega^{-1}_R\oplus \OOO_{C_\lambda}\cdot w_4 \Omega^{-1}_R,
\end{equation*}
where $w_3\Omega_R^{-1}= (\lambda w_1)^{-1} w_2\Omega_R^{-1}$. Whence, 
\begin{equation}
\label{eq:new-new-new:gr1Com}
 \gr^1_{C_\lambda}\upomega^*= (R)  \toplus \bigl(-1 +\textstyle \frac {m+3} 2 P^\sharp\bigr).
\end{equation}
For $\lambda\neq 0$ the germ $(X_\lambda, C_\lambda)$ is either of type \typec{kAD} or \typec{k2A_2}.  Comparing \eqref{eq:new-new-new:gr1Com} with \eqref{eq:new-new:k2A2} 
(resp. in view of \eqref{scase:kAD:new-1})
we see that $(X_\lambda, C_\lambda)$ is \typec{kAD} (resp. \typec{k2A_2}) if $m \ge 5$ (resp. $m=3$). \qqed
\end{proof}

\subsection{Case \typec{kAD}}
\label{case:kAD}
The configuration of singular points on
$(X,\, C)$ is the following: 
a type~\typec{IA} point $P$ of odd index $m\ge 3$ 
and a type~\typec{IA} point $Q$ of index $2$.
According to \cite[(A.3)]{Mori:flip}, \cite[(2.13)]{KM92}, and \cite{Mori-err} we can write 
\begin{equation*}
\begin{array}{cl}
(X,\, C,\, P)&=\left(\CC^3_{y_1,\, y_2,\, y_3},\, (y_1\axis),\, 0\right)/\mumu_{m}(1,(m+1)/2,-1),
\\[3pt] 
(X,\, C,\, Q)&=\bigl(\{\beta=0\}\subset \CC^4_{z_1,\dots,z_4},\, 
(z_1\axis),\, 0 \bigr) /\mumu_2(1,1,1,0),
\end{array}
\end{equation*}
where $\beta=\beta(z_1,\dots,z_4)$ is a semi-invariant with $\wt(\beta)\equiv 0\mod 2$.

For general divisor $S\in |-2K_X|$ we have $S\cap C=\{P\}$ and $S$ is given by $y_1=\xi(y_2,\, y_3)$ with $\xi \in (y_2,\, y_3)^2$ such that $\wt(\xi) \equiv 1 \mod m$. Thus 
\begin{gather}
S\simeq \CC_{y_2,\, y_3}^2/\mumu_m\left({\textstyle \frac {m+1}2},-1\right),\quad 
\omega_S=(\OOO_{S^\sharp,\, P^\sharp}\, \dd y_2\wedge \dd y_3)^{\mumu_m},\nonumber
\\
\label{eq:kAD:omegaCC}
\omega_S \otimes \CC_P = \CC \cdot y_2\, \dd y_2\wedge \dd y_3\oplus \CC\cdot y_3^{(m-1)/2}\, \dd y_2\wedge \dd y_3. 
\end{gather}
Then
\begin{equation}
\label{eq:new-new:gr0-b}
\gr_C^0\omega^*=\left(-1+{\textstyle \frac {m+1}2}P^\sharp+ Q^\sharp\right)\simeq \OOO_C(-1),
\end{equation}
where an $\ell$-free $\ell$-basis at $P$ and $Q$, respectively, can be written as follows: 
\begin{equation*}
\Omega_P^{-1}=\left(\dd y_1\wedge \dd y_2\wedge \dd y_3\right)^{-1},\qquad
\Omega_Q^{-1}=\left(\frac{\dd z_1\wedge \dd z_2\wedge \dd z_3}{\partial \beta/\partial z_4}\right)^{-1}
\end{equation*}
Hence, $H^0(C,\, \gr_C^0 \omega^*)=0$ and 
\begin{equation}
\label{eq:kAD:2:H0}
H^0\bigl(X,\, \OOO_X(-K_X)\bigr)=H^0\bigl(X,\, \I_C \totimes\OOO_X(-K_X)\bigr),
\end{equation}
where $\I_C$ is the defining ideal of $C$ in $X$.

As in \cite[(2.13.3)]{KM92} we distinguish two subcases:
\begin{align}
\label{scase:kAD:1}
\ell(Q)&\le 1,&i_Q(1)&=1, &\gr_C^1\OOO&\simeq \OOO\oplus \OOO(-1),
\\
\label{scase:kAD:2}
\ell(Q)&=2, &i_Q(1)&=2, &\gr_C^1\OOO&\simeq \OOO(-1)\oplus \OOO(-1).
\end{align}

\subsection{Subcase \eqref{scase:kAD:2}.}
\label{subcase:kAD:2}
It is treated similarly to~\ref{case:k3A}.
Since $\ell(Q)=2$, we have 
\begin{equation*}
\beta\equiv z_1^2z_4 \mod (z_2,\, z_3,\, z_4)^2. 
\end{equation*}
As in \cite[(2.13.4)]{KM92} we can arrange that
\begin{equation}
\label{eq:new-new:gr1-c}
\gr_C^1\omega^*=(0)\toplus \left(-1+{\textstyle \frac {m+3}2}P^\sharp\right),
\end{equation}
where 
\begin{equation}
\label{eq:kAD:basis}
\left(y_2\cdot \Omega_P^{-1},\, z_2\cdot \Omega_Q^{-1}\right),
\quad \left(y_3\cdot \Omega_P^{-1},\, z_3\cdot \Omega_Q^{-1}\right)
\end{equation}
form an $\ell$-free $\ell$-basis at $P$ and $Q$ for $(0)$ and $\left(-1+(m+3)/2P^\sharp\right)$, respectively, and 
\begin{equation*}
\beta \equiv z_1^2z_4+c_1z_3^2+c_2z_2z_3+c_3z_2^2 \mod \left(z_4,\, z_3^2,\, z_2z_3,\, z_2^2\right)\cdot \bigl(z_2,\, z_3,\, z_4\bigr)
\end{equation*}
for some $c_1,\, c_2,\, c_3 \in \CC$ such that 
$(c_1,\, c_2,\, c_3)\neq0$ if $m=3$ by the classification of 3-fold terminal singularities \cite[(6.1)]{Reid:YPG},
and $c_1 \neq0$ if $m \ge 5$ (see \cite[(2.12.6)]{KM92} and \cite[Remark~2]{Mori-err}).
The rest of the argument is the same as~\ref{case:k3A} (the type \typec{k3A}), except that 
we use \cite[(2.13.5)]{KM92} instead of \cite[(2.12.6)]{KM92}.

\begin{sremark}
\label{remark:kAD:2}
We note that the lowest power of $\mumu_2$-invariant variable $z_4$ that appear in $\beta$ (i.e. the axial multiplicity for $(X,\, P)$) remains the same for the defining equation of $D^\sharp$ under the elimination of variable of $\wt\equiv 1 \mod 2$. Thus the graph $\Delta(D,\, C)$ is as given as type \typec{kAD} in~\ref{theorem:ge}.
\end{sremark}

\begin{slemma}\label{new-lemma:1}
In the situation of \xref{case:kAD} with \eqref{scase:kAD:2}, let $(X_\lambda,C_\lambda)$ be the twisted extension \cite[1b.8.1]{Mori:flip} of the germ
\begin{equation*}
 (X_\lambda, Q)=\{\beta-\lambda z_4=0\}/\mumu_2 \supset (C_\lambda, Q)=(z_1\axis)/\mumu_2
\end{equation*}
by $u=(z_1z_2,\, z_1z_3)$. Then for $\lambda\neq 0$ the germ $(X_\lambda,C_\lambda)$ is of type \typec{k3A}.
\end{slemma}
\begin{proof}
When $0<|\lambda|\ll 1$,  a small  neighborhood  $X_\lambda\ni Q$ has two singular points on $C_\lambda$: 
a cyclic quotient at $Q$ and 
a Gorenstein point at $(\sqrt{\lambda},0,0,0)$.\qqed
\end{proof}

\subsection{Subcase \eqref{scase:kAD:1}}
\label{scase:kAD:new-1}
Note that in this case $m\ge 5$ (see \cite[(2.13.10)]{KM92} and \cite{Mori-err}).
Since $\ell(Q)\le 1$, we have 
\begin{equation*}
\beta\equiv z_4 \ \left(\text{resp. $z_1z_3+z_2^2$}\right)\mod \bigl(z_2^2,\, z_3,\, z_4\bigr)\bigl(z_2,\, z_3,\, z_4\bigr).
\end{equation*}
By \cite[(2.13.10)]{KM92}  we have
\begin{equation}\label{eq:new-new-new-gr1O}
\gr^1_C \OOO = \bigl(\textstyle\frac{m-1}{2}P^\sharp+Q^\sharp\bigr)\toplus \bigl(-1+P^\sharp+Q^\sharp\bigr)
\ \bigl(\text{resp. } \bigl(\frac{m-1}{2}P^\sharp\bigr)\toplus\bigl(-1+P^\sharp+Q^\sharp\bigr)\bigr).
\end{equation}
Tensoring it with \eqref{eq:new-new:gr0-b}  we obtain 
\begin{equation}
\label{eq:new-new:gr1C-omega}
\gr^1_C \omega^* = (1)\toplus \bigl(-1+\textstyle\frac{m+3}{2}P^\sharp\bigr)
\quad \bigl(\text{resp. } \bigl(Q^\sharp\bigr)\toplus\bigl(-1+\frac{m+3}{2}P^\sharp\bigr)\bigr),
\end{equation}
where $\bigl(y_2\Omega^{-1}_P,\, z_3\Omega^{-1}_Q\bigl)$ (resp. $\bigl(y_2\Omega^{-1}_P,\, z_4\Omega^{-1}_Q\bigl)$) is the $\ell$-free $\ell$-basis for the first $\ell$-summand of $\gr^1_C \omega^*$
and $\bigl(y_3\Omega^{-1}_P,\, z_2\Omega^{-1}_Q\bigl)$ for the second. 
Take the ideal $\J\subset \I$ as in \cite[(2.13.10-11)]{KM92}. Thus 
\begin{gather*}
(\I/\J)\totimes \omega^*=\bigl(-1+\textstyle\frac{m+3}2P^\sharp\bigr),\qquad H^0\bigl(X,\, \OOO_X(-K_X)\bigr)= H^0\bigl(\F^2(\omega^*,\, \J)\bigr),
\\
\J^\sharp =\left(y_3^2,\, y_2\right) \text{at $P$,}\qquad \J^\sharp =\left(z_2^2,\, z_3,\, z_4\right)\text{at $Q$.}
\end{gather*}
Then we have by \cite[(2.13.11)]{KM92}
\begin{equation*}
\gr^2(\omega^*,\, \J)=(0)\toplus \bigl(-1+ \textstyle\frac{m+5}2P^\sharp +Q^\sharp\bigr) 
\end{equation*}
where 
\begin{equation}
\label{eq:kAD:1:basis}
y_2\cdot \Omega_P^{-1},\quad y_3^2\cdot \Omega_P^{-1},\quad
z_3\cdot \Omega_Q^{-1},\quad 
z_2^2\cdot \Omega_Q^{-1} \ (\text{resp.} z_4\cdot\Omega_Q^{-1})
\end{equation}
form an $\ell$-free $\ell$-basis at $P$, $Q$. 
Thus we investigate
\begin{equation*}
H^0\bigl(\OOO_X(-K_X)\bigr) =H^0\bigl(F^2(\omega^*,\, \J)\bigr) \xlongrightarrow{\hspace*{17pt}} \gr_C^2(\omega^*,\, \J)
\end{equation*}
via the induced homomorphism
\begin{equation*}
H^0\bigl(\OOO_X(-K_X)\bigr) \xlongrightarrow{\hspace*{17pt}} \gr^2(\omega^*,\, \J)
\xlongrightarrow{\hspace*{17pt}} (\gr^2(\omega^*,\, \J))^\sharp \otimes \CC_{P^\sharp},
\end{equation*}
which, by $(y_1-\xi) \cdot (\gr^2(\omega^*,\, \J))^\sharp \otimes \CC_{P^\sharp} =0$, factors as
\begin{equation*}
H^0\bigl(\OOO_X(-K_X)\bigr) \xlongrightarrow{\hspace*{17pt}} \omega_S
\xlongrightarrow{\hspace*{17pt}} (\gr^2(\omega^*,\, \J))^\sharp \otimes \CC_{P^\sharp}
\end{equation*}
and further to
\begin{equation*}
\delta : H^0\bigl(\OOO_X(-K_X)\bigr) 
\xtwoheadrightarrow{\hspace*{10pt}} 
\omega_S\otimes \CC_P
\xlongrightarrow{\hspace*{17pt}}
(\gr^2(\omega^*,\, \J))^\sharp \otimes \CC_{P^\sharp}
\end{equation*}
since $\Omega_X^2 \subset \mathfrak m_{S,\, P}\cdot \omega_S$
as in the \typec{k3A} case.

The image of $\delta$ is generated by $y_2\cdot \Omega_P^{-1}$ if $m>5$, and 
by $y_2\cdot \Omega_P^{-1},\, y_3^2\cdot\Omega_P^{-1}$ if $m=5$
(see~\ref{eq:kAD:omegaCC} and~\ref{eq:kAD:1:basis}).
Hence if $\sigma_S$ is chosen general, the image $\bar \sigma$ of $\sigma$
in $\gr^2(\omega^*,\, \J)$ globally generates the direct summand $\OOO_C$ if $m>5$ and a general global section
of $\gr^2(\omega^*,\, \J) \simeq \OOO_C\oplus\OOO_C$ if $m=5$.
Hence,
\begin{equation*}
\bar \sigma\equiv 
\begin{cases}
(\lambda_Py_2+\mu_Py_3^2)\Omega_P^{-1}& \text{at $P$,}
\\
(\lambda_Qz_3+\mu_Qz_2^2)\Omega_Q^{-1}\ 
(\text{resp.}(\lambda_Qz_3+\mu_Qz_4)\Omega_Q^{-1})& \text{at $Q$,}
\end{cases}
\end{equation*}
where 
\begin{eqnarray*}
m\ge 7 &\Longrightarrow & \lambda_P(P) \lambda_Q(Q)\neq 0,
\\
m=5&\Longrightarrow &\text{$\lambda_P(P)$ and $\mu_P(P)$ are independent, and}
\\
& &\text{$\lambda_Q(Q)$ and $\mu_Q(Q)$ are independent.}
\end{eqnarray*}
These mean that 
the corresponding $D\in |{-}K_X|$ is smooth outside $P$ and $Q$ and
$D$ is Du Val at $P$ and $Q$ by \cite[(2.13.6)]{KM92}.
See Remark~\ref{remark:kAD:2} for further details.

\begin{slemma}\label{new-lemma:3}
In the situation of \xref{case:kAD} with \eqref{scase:kAD:1}, let $(X_\lambda,C_\lambda)$ be the twisted extension \cite[1b.8.1]{Mori:flip} of the germ
$ (X_\lambda, P)= (X,\, P) \supset (C_\lambda, P)=(C,\,  P)$
by $u=\bigl(y_1^{(m-1)/2}y_2+\lambda y_1y_3, \, y_1,y_3\bigr)$. Then for $\lambda\neq 0$ the germ $(X_\lambda,C_\lambda)$ is of type \typec{k2A_2}.
\end{slemma}
\begin{proof}
It is clear that $(X_\lambda,C_\lambda)$ is of type \typec{k2A_2} or \typec{kAD} with \eqref{scase:kAD:1} because $\ell(Q)\le 1$.
In either case, there is only one  non-zero section  $s$ (up to constant multiplication) of $\gr_C^1\OOO$  (cf. \eqref{eq:new-new:k2A2}). Since $s=y_1^{(m-1)/2}y_2\in \gr_C^1\OOO$ at $P$ (up to constant), we have an extension 
\[
s_\lambda = y_1^{(m-1)/2}y_2+\lambda y_1y_3= y_1\bigl(y_1^{(m-3)/2}y_2+\lambda y_3 \bigr)\in \gr_{C_\lambda}^1\OOO
\]
on $(C_\lambda, P)$ which generates $(P^\sharp)\subset \gr_{C_\lambda}^1\OOO$. In view of \eqref{eq:new-new-new-gr1O} the germ  $(X_\lambda,C_\lambda)$ is of type \typec{k2A_2} by \eqref{scase:kAD:new-1}. \qqed
\end{proof}

\begin{slemma}
In the situation of \xref{case:kAD} with \eqref{scase:kAD:1} and $\ell(Q)=1$, let $(X_\lambda,C_\lambda)$ be the twisted extension \cite[1b.8.1]{Mori:flip} of the germ
\begin{equation*}
 (X_\lambda, Q)=\{\beta-\lambda z_4=0\}/\mumu_2 \supset (C_\lambda, Q)=(z_1\axis)/\mumu_2
\end{equation*}
by $u=\bigl(z_1z_2,\, z_4\bigr)$. Then for $\lambda\neq 0$ the germ $(X_\lambda,C_\lambda)$ is of type \typec{kAD} and $\ell(Q)=0$.
\end{slemma}
Indeed, a global section $s$ of $\gr^1_C \OOO$ extends to $s_\lambda=y_4$ of $\gr^1_{C_\lambda} \OOO$ at $Q$, and $s_\lambda$ 
vanishes  at $P^\sharp$ to order $(m-1)/2>1$ by \eqref{scase:kAD:new-1}.
Thus $(X_\lambda, C_\lambda)$ is of type \typec{kAD}.

\subsection{Case \typec{k2A_2}}
\label{subcase:k2A}
This case comes from \cite[(2.13.1) and (2.13.9)]{KM92}.
The configuration of singular points on
$(X,\, C)$ is the following: 
a type~\typec{IA} point $P$ of odd index $m\ge 3$ 
and a type~\typec{IA} point $Q$ of index $2$.
According to \cite[(A.3)]{Mori:flip}, \cite[(2.13)]{KM92}, and \cite{Mori-err} we can write 
\begin{equation*}
\begin{array}{cl}
(X,\, C,\, P)&=\Bigl(\{\alpha=0\}\subset \CC^4_{y_1,\dots,y_4},\, (y_1\axis),\, 0\Bigr)/\mumu_{m}(1,a,-1,0),
\\[5pt] 
(X,\, C,\, Q)&=\Bigl(\{\beta=0\}\subset \CC^4_{z_1,\dots,z_4},\, 
(z_1\axis),\, 0 \Bigr) /\mumu_2(1,1,1,0),
\end{array}
\end{equation*}
where $a$ is an integer prime to $m$ such that $m/2<a<m$, and $\alpha$ and $\beta$
are invariants with
\begin{align*}
\alpha&=y_1y_3-\alpha_1(y_2,\, y_3,\, y_4),& \alpha_1 \in (y_2,\, y_3)^2+(y_4),
\\
\beta&=z_1z_3-\beta_1(z_2,\, z_3,\, z_4),& \beta_1 \in (z_2,\, z_3)^2+(z_4).
\end{align*}
Then
\begin{equation}
\gr_C^0\omega^*=\left(-1+aP^\sharp+ Q^\sharp\right)\simeq \OOO_C(-1),
\end{equation}
where an $\ell$-free $\ell$-basis at $P$ and $Q$, respectively, can be written as follows: 
\begin{equation*}
\Omega_P^{-1}=\left(\frac{\dd y_1\wedge \dd y_2\wedge \dd y_3}{\partial \alpha/\partial y_4}\right)^{-1}\qquad
\Omega_Q^{-1}=\left(\frac{\dd z_1\wedge \dd z_2\wedge \dd z_3}{\partial \beta/\partial z_4}\right)^{-1}
\end{equation*}
Hence, $H^0(C,\, \gr_C^0 \omega^*)=0$
and 
\begin{equation}
\label{eq:k2A:2:H0}
H^0\left(X,\, \OOO_X(-K_X)\right)=H^0\left(X,\, \I_C \totimes\OOO_X(-K_X)\right),
\end{equation}
where $\I_C$ is the defining ideal of $C$ in $X$.

As in the argument in \cite[(2.13.8--9)]{KM92}, we have
\begin{equation}
\label{eq:new-new:k2A2}
\gr_C^1\OOO=\LLL \toplus \gr_C^0\upomega,
\qquad
\gr_C^1\omega^*=\LLL\totimes (\gr_C^0 \omega^*)\toplus(0),
\end{equation}
where $\LLL$ is an $\ell$-invertible sheaf such that $\LLL=(P^\sharp+Q^\sharp)$ 
(resp. $(P^\sharp)$; $(Q^\sharp)$; $(0)$) 
if $y_4 \in \alpha$ and $z_4\in\beta$ 
(resp. $y_4\in\alpha$ and $z_4\not\in\beta$; $y_4\not\in\alpha$ and $z_4\in\beta$; $y_4\not\in\alpha$ and $z_4\not\in\beta$).
We also see that $y_3\Omega_P^{-1}$ (resp. $y_4\Omega_P^{-1}$) and $y_2\Omega_P^{-1}$
form an $\ell$-free $\ell$-basis for $\gr_C^1 \omega^*$ at $P$
if $y_4 \in \alpha$ (resp. $y_4\not\in \alpha$).

For a general divisor $S\in |-2K_X|$ we have $S\cap C=\{P\}$ and $S$ in $X$ is given by 
\begin{equation}\label{eq:new-gamma}
\gamma:=y_1^{2a-m}+y_2^2+y_3^{2m-2a}+\cdots=0 
\end{equation}
with $\wt(\gamma) \equiv 2a \mod m$. 
Let $\Omega$ be a generator of the dualizing sheaf $\upomega_{S^\sharp}$ of $S^\sharp$ at $P^\sharp$.
Then
\begin{gather}\label{eq:k2A:omegaCC1}
\omega_S=(\OOO_{S^\sharp,\, P^\sharp}\Omega)^{\mumu_m},\qquad \wt(\Omega) \equiv -a \mod m,
\\
\label{eq:k2A:omegaCC}
\omega_S \otimes \CC_P = \CC \cdot y_2\, \Omega \oplus \CC\cdot y_3^{m-a}\, \Omega.
\end{gather}

\begin{slemma}\label{lemma:new-zero}
The induced map 
\begin{equation}\label{eq:new-Omega}
\Omega_S^2 \xlongrightarrow{\hspace*{17pt}} \omega_S \otimes \CC_P
\end{equation}
is zero, where $\Omega_S^2$ is the sheaf of holomorphic $2$-forms on $S$. 
\end{slemma}

\begin{proof}
We have 
\begin{equation}\label{eq:newOmaga-Delta}
\Omega=\pm \frac{\dd y_1\wedge\dd y_2}{\varDelta_{3,4}}=\cdots=\pm \frac{\dd y_3\wedge\dd y_4}{\varDelta_{1,2}},
\qquad 
\varDelta_{i,j}:= 
\begin{vmatrix}
\frac{\partial \alpha}{\partial y_i}& \frac{\partial \alpha}{\partial y_j}
\\
\frac{\partial \gamma}{\partial y_i}& \frac{\partial \gamma}{\partial y_j}
\end{vmatrix}
\end{equation}
Note that $\wt(\varDelta_{i,j})\equiv 2a-\wt(y_i)-\wt(y_j) $ and $\wt(\Omega)\equiv -a \mod m$. Since $\upomega_S=\left(\OOO_{S^\sharp}\Omega\right)^{\mumu_m}$, it is sufficient to show that for any $\phi_1,\, \phi_2\in \CC\{y_1,\dots,y_4\}^{\mumu _m}$ 
the inclusion
\begin{equation}\label{eq:new:incl}
\dd \phi_1\wedge \dd \phi_2\in \mmm_{S,0}\cdot \left(\OOO_{S^\sharp}\Omega\right)^{\mumu_m}
\end{equation} 
holds. By \eqref{eq:newOmaga-Delta} the form $\dd \phi_1\wedge \dd \phi_2$ is a linear combination of the following 
\begin{equation*}
\frac{\partial(\phi_1,\, \phi_2)}{\partial y_i\partial y_j} \dd y_i\wedge \dd y_j= \frac{\partial(\phi_1,\, \phi_2)}{\partial y_i\partial y_j} \varDelta_{k,l}\Omega,
\quad 
\{i,j,k,l\}=\{1,\dots,4\}.
\end{equation*}
Denote 
\begin{equation*}
\Xi[i,j,k,l]:=\frac{\partial \phi_1}{\partial y_i} \cdot \frac{\partial \phi_2}{\partial y_j} \cdot \frac{\partial \alpha}{\partial y_k}\cdot \frac{\partial \gamma}{\partial y_l},
\quad \{i,j,k,l\}=\{1,\dots,4\}.
\end{equation*}
Since ${\partial(\phi_1,\, \phi_2)}/ (\partial y_i\partial y_j) \varDelta_{k,l}$ are linear combinations of $\Xi[i,j,k,l]$, it is sufficient to show that the following holds:
\begin{equation}\label{eq:new:Xi}
\left. \Xi[i,j,k,l]\ \right|_{S^\sharp}\in (\mmm_{S^\sharp})_{\wt=0}\cdot(\OOO_{S^\sharp})_{\wt=a}.
\end{equation}
First we note that 
\eqref{eq:new:Xi} holds if $\{1,\, 3 \}\subset \{i,\, j,\, k\}$. Indeed, let, for example, $i=1$, $j=3$.
Then 
\begin{equation*}
\frac{\partial \phi_1}{\partial y_1},\, \frac{\partial \phi_2}{\partial y_3} \in \mmm_{S^\sharp},
\qquad 
\wt \left( \frac{\partial \phi_1}{\partial y_1} \cdot \frac{\partial \phi_2}{\partial y_3} \right)\equiv0,
\qquad 
\wt \left( \frac{\partial \alpha}{\partial y_k} \cdot \frac{\partial \gamma}{\partial y_l} \right)\equiv a
\end{equation*}
(because $\{\wt(y_k),\, \wt(y_l)\} = \{0,\, a\}$). Thus, 
\begin{equation*}
(\partial \phi_1/\partial y_1)\cdot (\partial \phi_2/\partial y_3)\in (\mmm_{S^\sharp})_{\wt=0},\qquad 
(\partial \alpha/\partial y_k) \cdot (\partial \gamma/\partial y_l)\in(\OOO_{S^\sharp})_{\wt=a}.
\end{equation*}
Therefore, $l=3$ or $1$.

Let $l=3$. Then $\wt\left( \partial \gamma/\partial y_3\right) \equiv 2a+1$ and
\begin{equation*}
\wt \left( \partial \phi_1/\partial y_i \right), \ \wt\left( \partial \phi_2/\partial y_j \right), \ \wt\left( \partial \alpha /\partial y_k \right) \equiv -1,\ -a, \ 0 \end{equation*}
up to permutation of $i$, $j$, $k$. We claim that any product $\Pi$ of three monomials of weight $ -1, \, -a, \, 2a+1$  belongs to
$(\mmm_{S^\sharp})_{\wt=0} \cdot (\OOO_{S^\sharp})_{\wt=a}$. 
This is obvious if $ \Pi$ is divisible by $y_4$, $y_1 y_3 $, or $y_2$. 
So it is enough to consider the case $\Pi$ is a power of $y_1$ or $y_3$.
In the former case $\Pi$ is divisible by $y_1^{m-1} \cdot y_1^{m-a} \cdot y_1^{2a+1-m}= y_1^m \cdot y_1^a$, and in the latter $\Pi$ is divisible by $y_3 \cdot y_3^a \cdot y_3^{2m-2a-1}=y_3^m \cdot y_3^{m-a}$, which settles the claim.

Finally, let $l=1$. 
Similarly to the previous case, we show that any product $\Pi$ of monomials of weights $ -a,\, 1,\, 2a-1$ belongs to $(\mmm_{S^\sharp})_{\wt=0} \cdot (\OOO_{S^\sharp})_{\wt=a}$.
Again we can assume $\Pi$ is a power of $y_1$ or $y_3$. In the latter case, $\Pi$ is divisible by $y_3^a \cdot y_3^{m-1} \cdot y_3^{2m-2a+1}= y_3^{3m-a}=y_3^m \cdot y_3^{2m-a}$. In the former case, $\Pi$ similarly is divisible by $y_1^a$. 
By the equation \eqref{eq:new-gamma}, the monomial $y_1^{2a-m}$ belongs to $(y_2,\, y_3) \mmm_{S^\sharp}$, and we have $y_1^a \in (\mmm_{S^\sharp})_{\wt=0} \cdot (\OOO_{S^\sharp})_{\wt=a}$
by $a > 2a-m$. This concludes the proof of Lemma~\ref{lemma:new-zero}.\qqed
\end{proof}

By Lemma~\ref{lemma:new-zero} the homomorphism $\delta$ is factored as other cases:
\begin{equation}
\label{eq:k2A:delta}
\delta : H^0\bigl(\OOO_X(-K_X)\bigr) 
\xtwoheadrightarrow{\hspace*{10pt}} 
\omega_S\otimes \CC_P
\xlongrightarrow{\hspace*{17pt}}
(\gr_C^1 \omega^*)^\sharp \otimes \CC_{P^\sharp}.
\end{equation}
Thus we see that $\delta$ is surjective if and only if $a=m-1$ and $y_4 \in \alpha$.
If $\delta$ is not surjective, then  $y_2\Omega_P^{-1}$ generates its image which is   the second summand $(0)$ of $\gr_C^1 \omega^*$ and hence $\bar \sigma$ is a nowhere vanishing section.
The rest is the same as \cite[(2.13.9)]{KM92}.
This concludes the proof of Theorem~\ref{theorem:extending}.

\begin{proposition}\label{subcase:k2A+k1A:m}
Let $(X,\, \bar C)$ be an extremal curve germ whose central fiber $\bar C$ is reducible. 
Suppose that $\bar C$ contains a component $C$ of type \typec{k2A_2} and 
another component $C'$ of type \typec{k1A} meeting at a point $P$ of index $m>2$.
Assume further that $(X,\, \bar C)$ satisfies 
the condition~\xref{thm:main}\xref{thm:main:condition}.
Then a general member $D\in |{-}K_{X}|$ is Du Val in a neighborhood of $C \cup C'$.
\end{proposition}

\begin{proof}
The \typec{k2A_2} case~\ref{subcase:k2A} of Theorem~\ref{theorem:extending} shows that a general member $D (\supset C)$ of 
$|{-}K_{\bar{X}}|$ is Du Val in a neighborhood of $C$. 
If $D\not \supset C'$, then $D\cap C'=\{P\}$ because otherwise $D$ contains a Gorenstein point $P'$ of $X$ and so $D\cdot C'>1$
which contradicts $D\cdot C'=-K_X\cdot C'<1$ by \cite[(2.3.1)]{Mori:flip}, \cite[(3.1.1)]{MP:cb1}. 
Thus we may assume that $D\supset C'$.
We use the notation of~\ref{subcase:k2A}.
In view of \typec{IA} and \typec{IA^\vee} in \cite[A.3]{Mori:flip}, the fact that $D$ defined by $y_2+ \cdots=0$ contains $C'$ means that $(X,\, C')$ is of type \typec{IA}, $C'^\sharp$ is smooth at $P^\sharp$, and either $y_1$ or $y_3$ is a coordinate of $C'^\sharp$.

\begin{slemma}
$y_3$ is a coordinate of $C'^\sharp$, and hence we may assume that $C'^\sharp$ is the $y_3$-axis modulo a $\mumu_m$-equivariant change of coordinates.
\end{slemma}
\begin{proof}
Assume that $y_1$ is a coordinate of $C'^\sharp$. Then $\I_{C'}^\sharp$ is generated by 
\begin{equation*}
y_1^a \gamma_2+\delta_2,\qquad y_1^{m-1} \gamma_3 + \delta_3,\qquad y_1^m \gamma_4+\delta_4
\end{equation*}
with $\gamma_i \in \OOO_X$ and $\delta_i \in (y_2,\, y_3,\, y_4)^2\OOO_X^\sharp$. Thus 
\begin{equation*}
\I_C^\sharp + \I_{C'}^\sharp \subset (y_2,\, y_3,\, y_4,\, y_1^a)
\end{equation*}
and $\omega_X^\sharp \otimes \OOO_X^\sharp /(\I_C^\sharp + \I_{C'}^\sharp)$ contains a non-zero $\mumu_m$-invariant element $y_1^{m-a}\Omega_P$ by $m-a < a$. In view of the exact sequence
\begin{equation*}
0 \to \gr_{C \cup C'}^0 \omega \xlongrightarrow{\hspace*{17pt}} \gr_C^0 \omega \oplus \gr_{C'}^0 \omega \xlongrightarrow{\hspace*{17pt}} 
\left(\omega_X^\sharp \otimes \OOO_X^\sharp /(\I_C^\sharp + \I_{C'}^\sharp)\right)^{\mumu_m} \to 0,
\end{equation*}
we see that $H^1(\gr_{C \cup C'}^0 \omega) \neq0$. This implies that $(\bar{X},\, C\cup C')$ is a conic bundle germ and $C \cup C'$ is a whole fiber of the conic bundle \cite[Corollary~ 4.4.1]{MP:cb1}. However $(C+C' \cdot D) <2$, and this is impossible. \qqed
\end{proof}

From now on we assume that $C'^\sharp$ is the $y_3$-axis.
Hence $y_2,\, y_1$ (or $y_4$) form an $\ell$-free $\ell$-basis of $\gr_{C'}^1 \OOO$, and $y_2\Omega_P^{-1},\, y_1\Omega_P^{-1}$ (or $y_4\Omega_P^{-1}$) form an $\ell$-free $\ell$-basis of $\gr_{C'}^1 \omega^*$ at $P$. 
Furthermore, we see that
$\gr_{C'}^1(\omega^*)$ has a global section $\bar \sigma=(y_2+ \cdots)\Omega_P^{-1}$ induced by $\sigma$ defining $D$.
We also note $\gr_{C'}^0 \omega^*=((m-a)P^\sharp)$ since the weight $\wt'$ for $C'$ is $\wt' \equiv -\wt\mod m$.
According to the \typec{k2A_2} case~\ref{subcase:k2A} of Theorem~\ref{theorem:extending} the divisor $D$ is defined at $P$ by 
\begin{equation*}
y_2 \varPsi_1 + y_3^{m-a}\varPsi_2 =0
\end{equation*}
where $\varPsi_i$ are invariant functions by \eqref{eq:k2A:omegaCC1}-\eqref{eq:k2A:omegaCC}. We observe the surjections
\begin{equation*}
H^0\bigl(\OOO_X(-K_X)\bigr) 
\xtwoheadrightarrow{\hspace*{10pt}} 
\omega_S/\Omega_S^2
\xtwoheadrightarrow{\hspace*{10pt}} 
\omega_S\otimes \CC_P
\end{equation*}
by \eqref{prop:extension:cb:delta} and \eqref{prop:extension:cb:eq} for the former and by Lemma~\ref{lemma:new-zero} for the latter. In particular, $\Psi_2(P)\not=0$.
Since $C'=(y_3\axis)/\mumu_m$, we have $D\not\supset C'$. 
This contradicts our assumption.
Thus a general elephant $D$ of $(\bar{X},\, \bar{C})$ is proved to be Du Val in a neighborhood of $C \cup C'$. \qqed
\end{proof}

\section{Proof of the main theorem}
\begin{notation}
Let $(\bar{X},\, \bar{C})$ be an extremal curve germ with reducible central fiber $\bar{C}$ such that 
on each irreducible component $C_i$ of $\bar{C}$ there exist at most one point of index $>2$. Let $\{P_a\}_{a \in A}$ be the collection of such points.
For each $C_i$ without points of index $>2$, choose
one general point of $C_i$. Let $\{P_b\}_{b \in B}$ be
the collection of such points.
For each $i \in A \cup B$, let $S_i \in |-2K_{(X,\, P_i)}|$
be a general element on the germ $(X,\, P_i)$, and set $S=\sum_{i \in A \cup B} S_i$. Then
$S$ extends to an element $|-2K_X|$ by \cite[Thm. (7.3)]{Mori:flip}.
\end{notation}

\begin{zproof}[of Main Theorem {\xref{thm:main}}]
Take a general element $\sigma_i \in \OOO_{S_i}(-K_X)$, and 
\begin{equation*}
\sigma_S:=\sum_i \sigma_i \in \OOO_{S}(-K_X).
\end{equation*}
By~\ref{prop:extension}\ref{prop:extension:bir} or~\ref{prop:extension}\ref{prop:extension:cb:delta-bar}, the section
$\sigma_S \mod \Omega_S^2$ lifts to 
\begin{equation*}
s \in H^0\bigl(X,\, \OOO_X(-K_X)\bigr).
\end{equation*}
Let $C_{\mathrm{main}} \subset \bar C$ be the union of the irreducible component of type \typec{IC}, \typec{IIB}, \typec{kAD}, \typec{k3A}, or \typec{k2A_2}. 

By Theorem~\ref{theorem:extending}, the divisor $D:=\{s=0\} \supset C_{\mathrm{main}}$ is Du Val in a neighborhood of $C_{\mathrm{main}}$ and, 
for each irreducible component $C$ of $C_{\mathrm{main}}$, the graph $\Delta(D,\, C)$ is as described in~\ref{theorem:ge} by Theorem~\ref{theorem:extending}.
If $C_{\mathrm{main}}=\varnothing$, then we are done by Propositions~\ref{proposition:index2} and~\ref{proposition:super-easy-cases}. 
So we assume $C_{\mathrm{main}}\neq\varnothing$ and each irreducible component $C\subset \bar C$ intersects $C_{\mathrm{main}}$
(because singular points of $\bar C$ are non-Gorenstein on $\bar X$, see \cite[Cor. 1.15]{Mori:flip}, \cite[Lemma~4.4.2]{MP:cb1}). 
Then $D$ is normal and $C_{\mathrm{main}}$ is connected. If $C\not \subset D$, then 
\begin{equation*}
C \cap D \subset C_{\mathrm{main}} \cap D
\end{equation*}
and $D$ is Du Val in a neighborhood of $C$ as well because $C \cdot D < 1$ \cite[(0.4.11.1)]{Mori:flip} and $D$ is Cartier outside $C_{\mathrm{main}}$ \cite[Corollary~1.15]{Mori:flip}.

Suppose $C_{\mathrm{main}}$ contains a component of one of the types \typec{IC}, \typec{IIB}, \typec{kAD}, or \typec{k3A} and let $C \subset D$. 
Let $\upsilon: D\to D_0$ be the contraction of $C_{\mathrm{main}}$ and $P_0:=\upsilon(C_{\mathrm{main}})$. Then the point $P_0\in D_0$ is Du Val of type \type{D} or \type{E} (because $\upsilon$ is crepant and by Theorem~\ref{theorem:ge}).
Apply Lemma~\ref{lemma:surfaces} to $D_0$. We conclude that the surface $D_0$ has only Du Val singularities and the same is true for $D$.

If $C$ meets $C_{\mathrm{main}}$ at an index $2$ point, then $D\not\supset C$ by Lemma~\ref{lemma:index2point}.
The case where $C_{\mathrm{main}}$ consists of curves of type \typec{k2A_2} and $C$ intersects $C_{\mathrm{main}}$ at $P$ of index $m>2$ is treated in Proposition~\ref{subcase:k2A+k1A:m}. Theorem {\xref{thm:main}} is proved. \qqed
\end{zproof}

\begin{zproof}[of Corollary~\xref{cor:main}]
If $C_i\subset C$ is of type \typec{IIB}, then it contains a type \type{cAx/4}-point $P$ and 
all other components $C_j\subset C$ passing through $P$ are of types \typec{IIA} or \typec{II^\vee}. But according to 
\cite[Th. 6.4 \& 9.4]{Mori:flip} $P$ is the only non-Gorenstein point on $C_j$. Since $X$ is not Gorenstein at any intersection point $C_i\cap C_j$ by \cite[Cor. 1.15]{Mori:flip} and \cite[Lemma 4.4.2]{MP:cb1},
all the components $C_k\subset C$ must pass through $P$. This proves~\ref{cor:main1}.

From now on we may assume that $C=C_i\cup C_j$. We also may assume that $C_j$ is not of type  \typec{k2A_{n,m}}, $n,\, m\ge 3$ (otherwise there is nothing to prove). Thus $(X,\, C)$ satisfies the condition \ref{thm:main}\ref{thm:main:condition}.
Let $f: (X,\, C)\to (Z,\, o)$ be the corresponding contraction.
Consider a general member $D\in |{-}K_X|$ and the Stein factorization
\begin{equation}
\label{eq:Stein}
f_D: D\supset C \xlongrightarrow{\hspace*{10pt} f'\hspace*{10pt}} D_Z\ni o_Z \xlongrightarrow{\hspace*{17pt}} f(D)\ni o. 
\end{equation}
By Theorem ~\xref{thm:main} the surface $D$ has only Du Val singularities.
The contraction $f': D\to D_Z$ is crepant and point $D_Z\ni o_Z$ is Du Val.
Now we note that the germs $(D,\, C_i)$ and $(D,\, C_j)$ are as described by Theorem~\ref{theorem:ge}.
Thus the whole configuration $\Delta(D,\, C)$ is one of the Dynkin diagrams \type{A}, \type{D} or \type{E}. 
In particular, $\Delta(D,\, C)$ has no vertices of valency $\ge 4$ and at most one vertex, say $v$, of valency 3.
On the other hand, by Theorem~\ref{theorem:extending} the configuration $\Delta(D,\, C)$ is obtained by ``gluing'' the configurations $\Delta(D,\, C_i)$ described in Theorem~\ref{theorem:ge} along one connected component of white subgraph.
Since the whole configuration $\Delta(D,\, C)$ is Du Val, at most one component of $C$ is of type \typec{IC}, \typec{IIB}, \typec{kAD} or \typec{k3A}.

For~\ref{cor:main2} it is enough  to note that all the singularities along $C_i$ are of type \type{cA} and so  extremal germs \type{cD/m}, \type{cAx/2}, \type{cE/2}, \typec{IIA}, and \typec{II^\vee} do not occur as a component of $(X,\, C)$.
Similar argument is applied in~\ref{cor:main3} but in this case singularities \type{cD/2} and \type{cAx/2} are allowed \cite{Mori-err}. 
It remains to prove  \ref{cor:main4}.

\subsection{}
\label{extend-deform-Cj-to-Ci}
Let $T$ be a point of $C_j$. It is easy to observe that a twisted extension $(X_{j,\lambda}, C_{j,\lambda})$ \cite[1b.8.1]{Mori:flip} 
of the germ $(X_{j,\lambda},T) \supset (C_{j,\lambda}, T)$ by $u$ in a neighborhood of  
$C_{j,\lambda}$ can naturally contain a neighborhood of $C_i$ if 
\begin{enumerate}
\renewcommand\labelenumi{\rm \theequation\alph{enumi}}
\renewcommand\theenumi{\rm \theequation\alph{enumi}}
\item\label{new:def:case(1)}
$T \not\in C_i$ or 
\item\label{new:def:case(2)}
$(X_{j,\lambda},T) \supset (C_{j,\lambda}, T)$ is a trivial deformation, i.e. $X_{j,\lambda}=X_j$ and $C_{j,\lambda}=C_j$.
\end{enumerate}

We can make a successive deformation of $(X,C)$ in a neighborhood of $C_j$ which is trivial on a neighborhood of $C_i$ in such a way that $(X,C_j)$ deforms 
\begin{enumerate}
\renewcommand\labelenumi{\rm (def.\arabic{enumi})}
\renewcommand\theenumi{\rm (def.\arabic{enumi})}
\item
\label{new:def:case(a)}
from \typec{kAD} of the case \eqref{scase:kAD:2} to \typec{k3A}, as treated in  Lemma~\ref{new-lemma:1}; 
\item
\label{new:def:case(b)}
from \typec{k3A} to \typec{k2A_2} or \typec{kAD} of the case \eqref{scase:kAD:1}, as treated in  Lemma~\ref{new-lemma:2}; 
\item
\label{new:def:case(c)}
from \typec{kAD} of the case \eqref{scase:kAD:1} to \typec{k2A_2}, as treated in Lemma~\ref{new-lemma:3},
\end{enumerate}
ultimately to $(X,\,  C_j)$ of type \typec{k2A_2}.
Indeed, \ref{new:def:case(1)} works with $T$ as the type \typec{III} point $\not\in C_i$  for deformation \ref{new:def:case(b)};
\ref{new:def:case(2)}  with $T$ as the index $m$ point for deformation \ref{new:def:case(c)}.

\subsection{}
For deformation \ref{new:def:case(a)}, we need to take as $T$ the index $2$ point $Q$ of $\ell(Q)=2$.
Suppose $Q \in C_i$, in which case the divisor $D$  cannot be Du Val because $\Delta(D,C_j)$ of type \type{D_k} with $k\ge 8$ and $\Delta(D,C_i)$ of type \type{A_q} with $q \ge 4$ are connected at the index 2 point $Q$. So $Q\not\in C_i$ and \ref{new:def:case(1)} applies, and we are left with the case $C_j$ is of type \typec{k2A_2}.
It remains to disprove the case where  both components of $C$ are of type \typec{k2A_2}.
This follows from Lemma~\ref{lemma:k2A:chain2} below. \qqed
\end{zproof}

\begin{lemma}\label{k2A-H}
Let $(X,\, C)$ be an extremal curve germ such that any component $C_i\subset C$ is of type \typec{k2A}.
Assume that a general member $D\in |{-}K_X|$ is Du Val. 
Then a general hyperplane section $H\in |\OOO_X|_C$ passing through $C$ has only cyclic quotient singularities and the pair $(H,\, C)$ is log canonical and purely log terminal outside $\Sing(C)$. 
\end{lemma}

\begin{sremark}
If in the conditions of Lemma~\ref {k2A-H} the germ $(X,\, C)$ is a $\QQ$-conic bundle, then its base surface is smooth.
Indeed,
if the base surface is singular, then by \cite[Theorem~1.3]{MP:cb2} each component $C_i\subset C$ must be locally imprimitive. 
\end{sremark}

\begin{proof}
Let $f: (X,\, C)\to (Z,\, o)$ be the corresponding contraction. Note that $D\supset C$. Consider the Stein factorization \eqref{eq:Stein}.
It is easy to see that the configuration $\Delta(D,\, C)$ is a linear chain. Therefore, $D_Z\ni o_Z$ is a (Du Val) singularity of type \type{A}. Then the arguments in the proof of \cite[Proposition~2.6]{MP:IA} work and show that the pair $(X,D+H)$ is log canonical. Since $D\supset C$, the pair $(X,H)$ is purely log terminal by Bertini's theorem. Hence $H$ is normal and 
by the inversion of adjunction the pair $(H,\, C)$ is log canonical. Moreover, $C=D\cap H$ and $K_H+C$ is a Cartier divisor on $H$. By the classification of surface log canonical pairs 
\cite[Theorem~4.15]{KM:book} 
the singularities of $H$ are cyclic quotients and the pair $(H,\, C)$ is purely log terminal outside $\Sing(C)$. \qqed
\end{proof}

\begin{lemma}\label{lemma:k2A:chain2}
Let $(X,\, C)$ be an extremal curve germ, where  $C$ is reducible and has exactly two components. 
Then both components  cannot be of  type \typec{k2A_2}.
\end{lemma}

\begin{proof}
Assume the contrary. 
The computation below are very similar to that in \cite[Prop.~2.6]{Mori:ss}.
Let $C=C_1\cup C_2$, let $C_1\cap C_2=\{P_0\}$, and let $P_i\in C_i$ for $i=1$, $2$ be the non-Gorenstein point other than $P_0$.
Let $H\in |\OOO_X|_C$ be a general hyperplane section  passing through $C$.
According to Lemma~\ref{k2A-H} the surface $H$
has only cyclic quotient singularities. 
Consider the minimal resolution $\mu: \tilde H\to H$ and write 
\begin{equation}
\label{eq:ss:KC}
K_{\tilde H}=\mu^*K_H-\Theta,
\end{equation}
where $\Theta$ is an effective $\QQ$-divisor with support in the exceptional locus (codiscrepancy divisor).
Let $\tilde C_i$ be the proper transform of $C_i$. Then $K_{\tilde H}\cdot \tilde C_i<K_H\cdot C_i<0$.
Hence $\tilde C_i$ is a $(-1)$-curve on $\tilde H$. Moreover, 
\begin{equation*}
\Theta \cdot \tilde C_i=1+K_H\cdot C_i<1.
\end{equation*}
Since $H$ is a Cartier divisor on $X$ such that $(X,H)$ is purely log terminal, the singularities of $H$ are of type \type{T} \cite{KSh88}. 
Hence,
\begin{equation*}
(H\ni P_i)  \simeq \left(\CC/\mumu_ {m_i^2 p_i}(1, m_i p_i a_i-1) \ni 0 \right)
\end{equation*}
for some positive $m_i$, $p_i$, $a_i$
such that $a_i<m_i$ and $\gcd (m_i,a_i)=1$. Here $m_i$ is the index of $P_i$. Write $\Theta=\Theta_0+\Theta_1+\Theta_2$ so that $\Supp(\Theta_i)=\mu^{-1}(P_i)$. Computations with weighted blowups \cite[10.1-10.3]{KM92} show that the coefficients of $\Theta_i$ in the ends of the chain $\Supp(\Theta_i)$ are equal to $(m_i-a_i)/m_i$ and $a_i/m_i$.
Since $\tilde C_1$ and $\tilde C_2$ meet different ends of the chain $\Supp(\Theta_0)$, up to permutations $P_1$ and $P_2$ and changing generators of $\mumu_ {m_i^2 p_i}$ we have 
\begin{equation*}
\Theta_1\cdot\tilde C_1=\textstyle \frac{a_1}{m_1},\quad \Theta_0\cdot \tilde C_1= \frac{m_0-a_0}{m_0},\quad 
\Theta_0\cdot \tilde C_2= \frac{a_0}{m_0},\quad \Theta_2\cdot \tilde C_2= \frac{m_2-a_2}{m_2}.
\end{equation*}
Denote
\begin{equation}\label{eq:c:delta}
\delta_1:= a_0m_1 -a_1 m_0,\qquad \delta_2:= a_2m_0 -a_0 m_2.
\end{equation}
Then by \eqref{eq:ss:KC}
\begin{equation*}\textstyle 
-K_H\cdot C_1=\frac{a_0}{m_0} -\frac{a_1}{m_1}= \frac{\delta_1} {m_0m_1}>0,\quad
-K_H\cdot C_2=\frac{a_2}{m_2} -\frac{a_0}{m_0}= \frac{\delta_2} {m_0m_2}>0.
\end{equation*}
Further, put
\begin{equation}\label{eq:Delta1-2}
\varDelta_i:=m_0^2 p_0+m_i^2 p_i -m_0p_0m_ip_i \delta_i,\quad i=1,\ 2.
\end{equation}
Then 
\begin{equation*}\textstyle 
C_i^2 =\frac{-\varDelta_i}{m_0^2 p_0m_i^2 p_i},\qquad 
C_1\cdot C_2= \frac1{m_0^2 p_0}.
\end{equation*}
Since the configuration is contractible, we have $\varDelta_i>0$ and
\begin{equation}\label{eq:det-Delta}
\varDelta_1\varDelta_2-m_1^2 p_1 m_2^2 p_2\ge 0.
\end{equation}
Assume $m_0 > 2$. Since $C_i$ is of type \typec{k2A_2}, $m_1 =m_2= 2$. Then
$a_1=a_2=1$ and \eqref{eq:c:delta} implies
\begin{equation*}
\delta_1= 2a_0 - m_0>0,\qquad \delta_2= m_0 -2a_0>0,
\end{equation*}
which is impossible. Hence, $m_0=2$. Then $a_0=1$ and 
\eqref{eq:Delta1-2} can be written as follows
\begin{equation*}
\varDelta_i=m_i^2 p_i -2p_0(m_ip_i \delta_i-2)>0,\quad i=1,\ 2,
\end{equation*}
where $m_ip_i \delta_i-2>0$. Then \eqref{eq:det-Delta} reads
\begin{equation*}
2p_0(m_1p_1 \delta_1-2)(m_2p_2 \delta_2-2) \ge (m_1p_1 \delta_1-2)m_2^2 p_2+(m_2p_2 \delta_2-2)m_1^2 p_1.
\end{equation*}
Combining these inequalities we obtain
\begin{equation*}\textstyle 
\frac{m_1^2 p_1}{m_1p_1 \delta_1-2}> 2p_0 \ge 
\frac{m_2^2 p_2}{m_2p_2 \delta_2-2}
+ \frac{m_1^2 p_1}{m_1p_1 \delta_1-2}
>  \frac{m_1^2 p_1}{m_1p_1 \delta_1-2}.
\end{equation*}
The contradiction proves the lemma. \qqed
\end{proof}

\end{fulltext}

\appendix

\section{Examples}
\renewcommand{\thesubsection}{\Alph{section}.\arabic{subsection}}
\subsection{}
As was explained in \cite[Sect. 6.6]{MP-1p}, to construct an example of an extremal curve germ whose general member of $|\OOO_X|$ is normal, it is sufficient to construct a normal surface germ $(H,\, C)$ along a proper connected curve $C$ satisfying the following conditions:
\begin{enumerate}
\renewcommand\labelenumi{\rm (\alph{enumi})}
\renewcommand\theenumi{\rm (\alph{enumi})}
\item 
there exists a contraction $f_H: (H,\, C)\to (T, o)$ such that $-K_H$ is relatively ample and $f_H^{-1}(o)_{\red}=C$,
\item
for every each singular point $P_i\in H$ there exists a threefold terminal singularity $(X_i,\, P_i)$ and an embedding $(H,\, P_i)\subset (X_i,\, P_i)$ such that $H$ is a Cartier divisor on $X_i$ (at $P_i$). 
\end{enumerate}
Then there exists an extremal curve germ $(X,\, C)$ and an embedding $(H,\, C)\subset (X,\, C)$ such that $H$ is a member of $|\OOO_X|$ on $X$. 
Moreover, near each point $P_i\in H \subset X$ there is an isomorphism $(X_i,\, P_i)\simeq (X,\, P_i)$ inducing the embedding $(H,\, P_i)\subset (X_i,\, P_i)$:
\begin{equation*}
\xymatrix@R=7pt{
(X_i,\, P_i)\ar@{}[]!<2pt,-2pt>;[d]!<2pt,-2pt> |-*[@]{\supset} \ar@{}[]!<2pt,-2pt>;[r]!<2pt,-2pt> |-*[@]{\simeq} & (X,\, P_i)\ar@{}[]!<2pt,-2pt>;[d]!<2pt,-2pt> |-*[@]{\supset}
\\
(H,\, P_i)\ar@{}[]!<2pt,-2pt>;[r]!<2pt,-2pt> |-*[@]{=} &(H,\, P_i)
} 
\end{equation*}
Then for any component $C_j\subset C$, types of points $P_i\in (X,\, C_j)$ can be seen locally near each $P_i$.

\subsection{}
To construct examples of an extremal curve germs with reducible central curve one can start with one of the configurations described in the Appendix of \cite{KM92} or \cite{MP:ICIIB} and add extra black vertices. We need to check only 
the ampleness of the anti-canonical divisor. Below we give a few examples. Of course, they do not exhaust all the possibilities.
Below in the dual graphs of minimal resolutions we use the notation~\ref{not:graphs}.
A number attached to a white vertex denotes the negative of the self-intersection number.
We may omit $2$ if the self-intersection equals $-2$.

\begin{example}
Starting with an exceptional flipping extremal curve germ $(X,\, C_1)$ of type \typec{IC} as in \cite[A.3.2.1]{KM92} we can construct the following configurations:
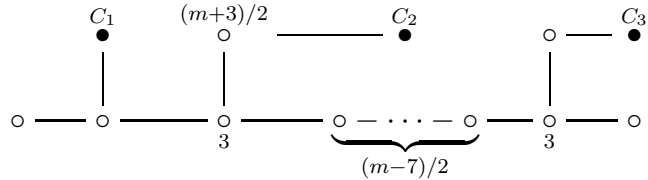
\begin{equation*}
\vcenter{
\xymatrix@R=4pt@C=17pt{
&\overset{C_1}{\bullet}\ar@{-}[d]&\overset{(m+3)/2}\circ\ar@{-}[d]&\overset{C_2}{\bullet}\ar@{-}[l]
&\circ\ar@{-}[d]&\overset{C_3}{\bullet}\ar@{-}[l]
\\
\circ\ar@{-}[r] &\circ\ar@{-}[r] &\underset3\circ\ar@{-}[r] &{\underbrace{\circ -\cdots - \circ}_{(m-7)/2}}&\underset3\circ\ar@{-}[l]\ar@{-}[r]&\circ
} }
\leqno{\mathrm{}}
\end{equation*}
Computing the intersection numbers on the minimal resolution of $H$ we obtain $-K_H\cdot C_1=-K_H\cdot C_2=1/m$ and $-K_H\cdot C_3=(m-1)/2m$.
Each $\bullet$ corresponds to a $(-1)$-curve.
Hence, $-K_H$ is ample and our construction gives an example of a divisorial contraction of type \typec{IC}+\typec{k1A}+\typec{k1A}. 
(see Corollary~\ref{cor:main}\ref{cor:main2}).
\end{example}

\begin{example}
Starting with an exceptional divisorial extremal curve germ $(X,\, C_1)$ of type \typec{IIB} as in \cite[1.2.2]{MP:ICIIB} we can construct the following configurations:
\begin{equation*}
\vcenter{
\xymatrix@R=4pt@C=17pt{
\overset{3}\circ\ar@{-}[r]&\overset{4}{\circ}\ar@{-}[d]\ar@{-}[r]&{\circ}\ar@{-}[d]\ar@{-}[r]&\circ\ar@{-}[r]&\circ
\\
\underset{C_2}{\bullet}\ar@{-}[r]&\underset{3}{\circ}&\circ\ar@{-}[r]&\underset{C_1}{\bullet}
}}
\end{equation*}
Computing the intersection numbers on the minimal resolution of $H$ we obtain $-K_H\cdot C_1=-K_H\cdot C_2=1/4$.
Hence, $-K_H$ is ample and the germ $(X,\, C_2)$ is locally primitive.
Therefore, our construction gives an example of a divisorial contraction of type \typec{IIB}+\typec{IIA} (see Corollary~\ref{cor:main}\ref{cor:main1}).
\end{example}

\begin{example}
Starting with an exceptional divisorial extremal curve germ $(X,\, C_1)$ of type \type{cD/3} as in \cite[4.5.2.1]{MP:IA} we can construct the following configurations:
\begin{equation*}
\xymatrix@R=3pt{
\circ\ar@{-}[r]&\overset{C_1}{\bullet}\ar@{-}[r]&\overset{3}{\circ}\ar@{-}[r]&\circ\ar@{-}[r]&\overset{3}{\circ}\ar@{-}[r]&\overset{C_2}{\bullet}
\\
&&&\underset{3}{\circ}\ar@{-}[r]\ar@{-}[u]&\underset{C_3}{\bullet}
} 
\end{equation*}
So, this is an example of a $\QQ$-conic bundle contraction of type \type{cD/3}+\type{cD/3}+\type{cD/3}. 
\end{example}

\section{$\QQ$-conic bundles with irreducible central fiber}
\label{appB}
For convenience of references, we provide below  the classification of  irreducible
$\QQ$-conic bundle extremal curve germs and their general anticanonical members.

\subsection*{Notation}
Let $(X,C)$ be a $\QQ$-conic bundle extremal curve germ with irreducible $C$, let $f:X\to Z$ be the corresponding contraction, and  let $D\in |-K_X|$ be the general member. If $D\supset C$, then 
$\Delta(D,C)$ denotes the dual graph of the minimal resolution $\mu: \tilde D\to D$.
For this graph, we use the standard notation: each vertex $\circ$ corresponds to
a prime $\mu$-exceptional divisor (which is a $(-2)$-curve on $\tilde D$) and 
vertex $\bullet$ corresponds to the proper transform of $C$ (which is a $(-1)$-curve).

\subsection*{The cases with singular $Z$ \cite{MP:cb1}}

\subsubsection*{Case~\typec{T}, see \cite[1.2.1]{MP:cb1}}
$D$ does not contain $C$ and it is a disjoint union $D=D_1+D_2$, where $D_1\simeq D_2$ is a 
singularity of type~\typet{A}{m-1}, $m\ge 2$.

\subsubsection*{Case~\typec{k2A}, see \cite[1.2.2 and Theorem~11.1]{MP:cb1}}
$D\supset C$, it is given by $y_1=0$, and $\Delta(D,C)$ is as follows:
\[
\begin{tikzpicture}
 \def\sizec{0.2em}
\def\sl{0.7}
\coordinate (Q2) at (1,1*\sl);
\coordinate (P1) at (0,0);
\coordinate (P2) at (1,0);
\coordinate (P3) at (2,0);
\coordinate (P4) at (3,0);
\coordinate (P5) at (4,0);
\coordinate (P6) at (5,0);
\coordinate (P7) at (6,0);
\coordinate (P8) at (7,0);
\coordinate (P9) at (8,0);
\coordinate (P10) at (9,0);
\coordinate (P11) at (10,0);
\coordinate (Q12) at (11,0);

\draw (P2) -- (P3);
\draw [loosely dotted] (P3) -- (P4);
\draw (P4) -- (P5);
\draw (P5) -- (P6);
\draw (P6) -- (P7);
\draw (P7) -- (P8);
\draw [loosely dotted](P8) -- (P9);
\draw (P9) -- (P10);

\path[fill=white,draw=black] (P2) circle (\sizec)node [below, yshift=0, xshift=0] {};
\path[fill=white,draw=black] (P5) circle (\sizec)node [below, yshift=0, xshift=0] {};
\path[fill=black,draw=black] (P6) circle (\sizec)node [below, yshift=0, xshift=0] {};

\path[fill=white,draw=black] (P7) circle (\sizec)node [below, yshift=0, xshift=0] {};

\path[fill=white,draw=black] (P10) circle (\sizec)node [below, yshift=0, xshift=0] {};

\draw[thick,black,decorate,decoration={brace,amplitude=10pt,mirror, 
raise=5pt},xshift=0.4pt,yshift=-1.4pt](P2) -- (P5) node[black,midway,yshift=-23] {\footnotesize 
$m-1$};
\draw [thick, black,decorate,decoration={brace,amplitude=10pt,mirror, 
raise=5pt},xshift=0.4pt,yshift=-1.4pt](P7) -- (P10) node[black,midway,yshift=-23] {\footnotesize 
$m-1$};
\end{tikzpicture}
\]
$m$ is odd $\ge 3$.
\subsubsection*{Case~\typec{IE^\vee}, see \cite[1.2.3]{MP:cb1}}
The intersection 
$D\cap C$ is a single point $P$ and $(D, P)$ is of type~\typet{A}{7}.

\subsubsection*{Case~\typec{ID^\vee}, see \cite[1.2.4]{MP:cb1}}
$D\cap C=\{P\}$ and $(D, P)$ is of type~\typet{A}{3} or 
\typet{D}{k}, $k\ge 4$.

\subsubsection*{Case~\typec{IA^\vee}, see \cite[1.2.5]{MP:cb1}}
$D\cap C=\{P\}$ and $(D, P)$ is of type~\typet{A}{3}.

\subsubsection*{Case~\typec{II^\vee}, see \cite[1.2.6]{MP:cb1}}
$D\cap C=\{P\}$ and $(D, P)$ is of type~\typet{D}{2k+1}, $k\ge 2$.

\subsection*{The cases with  smooth $Z$ \cite{MP:cb1}, \cite{MP:cb3}, 
\cite{MP:IA}, \cite{MP:ICIIB}}
\subsubsection*{Gorenstein case}
$|-K_X|$ is base point free and $(D, P)$ is smooth.

\subsubsection*{Case~\typec{cAx/2}, see \cite[\S~12]{MP:cb1}, 
\cite[\S~7]{MP:IA}}
$D\cap C=\{P\}$ and $(D, P)$ is of type~\typet{D}{k}, $k\ge 4$.

\subsubsection*{Case~\typec{cD/2}, see \cite[\S~12]{MP:cb1}, 
\cite[\S~7]{MP:IA}}
$D\cap C=\{P\}$ and $(D, P)$ is of type~\typet{D}{2k}, $k\ge 3$.

\subsubsection*{Case~\typec{cE/2}, see \cite[\S~12]{MP:cb1}, \cite[\S~7]{MP:IA}}
$D\cap C=\{P\}$ and $(D, P)$ is of type~\typet{E}{7}.

\subsubsection*{Case~\typec{k1A}, see \cite[Theorem~8.6(iii)]{MP:cb1}}
$D\cap C=\{P\}$ and $(D, P)$ is of type~\typet{A}{k}, $k\ge 1$ 

\subsubsection*{Case~\typec{cD/3}, see \cite[Theorem~8.6(iii)]{MP:cb1}}
$D\cap C=\{P\}$ and $(D, P)$ is of type~\typet{E}{6} 

\subsubsection*{Case~\typec{IC}, see \cite[Theorem~1.3]{MP:cb3} and \cite[Theorem~1.1]{MP:ICIIB}}
$D\supset C$, the unique non-Gorenstein point has index $5$, and $\Delta(D,C)$ is as follows:
\[
\begin{tikzpicture}[square/.style={regular polygon,regular polygon sides=4}]
\def\sizec{0.2em}
\def\sl{0.7}
\coordinate (1) at (0,0);
\coordinate (2) at (1,0);
\coordinate (3) at (2,0);
\coordinate (4) at (3,0);
\coordinate (5) at (4,0);
\coordinate (6) at (1,1*\sl);
\coordinate (7) at (2,1*\sl);
\coordinate (8) at (3,1*\sl);
\coordinate (9) at (2,2*\sl);

	\path (2) edge (3);
	\path (3) edge (4);
	\path (4) edge (5);
	\path (8) edge (4);
\path[fill=white,draw=black] (2) circle (\sizec) ;
\path[fill=white,draw=black] (3) circle (\sizec);
\path[fill=white,draw=black] (4) circle (\sizec) ;
\path[fill=black,draw=black] (5) circle (\sizec) ;
\path[fill=white, draw=black] (8) circle (\sizec);
\end{tikzpicture}
\]

\subsubsection*{Case~\typec{IIA}, see \cite{MP:IIA-1}}
$D\cap C=\{P\}$ and $(D, P)$ is of type~\typet{D}{2k+1}, $k\ge 2$.

\subsubsection*{Case~\typec{IIB}, see \cite[Theorem~1.3]{MP:cb3}}
$D\supset C$ and $\Delta(D,C)$ is as follows:
\[
\begin{tikzpicture}[square/.style={regular polygon,regular polygon sides=4}]
\def\sizec{0.2em}
\def\sl{0.7}
\coordinate (1) at (0,0);
\coordinate (2) at (1,0);
\coordinate (3) at (2,0);
\coordinate (4) at (3,0);
\coordinate (5) at (4,0);
\coordinate (6) at (1,1*\sl);
\coordinate (7) at (2,1*\sl);
\coordinate (8) at (3,1*\sl);
\coordinate (9) at (2,2*\sl);
	\path (1) edge (2);
	\path (2) edge (3);
	\path (3) edge (4);
	\path (4) edge (5);
	\path (7) edge (3);
	
	\path[fill=white,draw=black] (1) circle (\sizec) ;
\path[fill=white,draw=black] (2) circle (\sizec) ;
\path[fill=white,draw=black] (3) circle (\sizec);
\path[fill=white,draw=black] (4) circle (\sizec) ;
\path[fill=black,draw=black] (5) circle (\sizec) ;
\path[fill=white, draw=black] (7) circle (\sizec);
\end{tikzpicture}
\]

\subsubsection*{Case~\typec{kAD}, see \cite[Theorem~1.3]{MP:cb3}}
$D\supset C$ and $\Delta(D,C)$ is as follows:
\[
\begin{tikzpicture}
 \def\sizec{0.2em}
\def\sl{0.7}
\coordinate (Q2) at (9,1*\sl);
\coordinate (P1) at (0,0);
\coordinate (P2) at (1,0);
\coordinate (P3) at (2,0);
\coordinate (P4) at (3,0);
\coordinate (P5) at (4,0);
\coordinate (P6) at (5,0);
\coordinate (P7) at (6,0);
\coordinate (P8) at (7,0);
\coordinate (P9) at (8,0);
\coordinate (P10) at (9,0);
\coordinate (P11) at (10,0);
\coordinate (Q12) at (11,0);

\draw (P2) -- (P3);
\draw [loosely dotted] (P3) -- (P4);
\draw (P4) -- (P5);
\draw (P5) -- (P6);
\draw (P6) -- (P7);
\draw (P7) -- (P8);
\draw [loosely dotted](P8) -- (P9);
\draw (P9) -- (P10);
\draw (P10) -- (P11);
\draw (Q2) -- (P10);

\path[fill=white,draw=black] (P2) circle (\sizec)node [below, yshift=0, xshift=0] {};
\path[fill=white,draw=black] (P5) circle (\sizec)node [below, yshift=0, xshift=0] {};
\path[fill=black,draw=black] (P6) circle (\sizec)node [below, yshift=0, xshift=0] {};

\path[fill=white,draw=black] (P7) circle (\sizec)node [below, yshift=0, xshift=0] {};

\path[fill=white,draw=black] (P10) circle (\sizec)node [below, yshift=0, xshift=0] {};
\path[fill=white,draw=black] (P11) circle (\sizec)node [below, yshift=0, xshift=0] {};
\path[fill=white,draw=black] (Q2) circle (\sizec)node [below, yshift=0, xshift=0] {};

\draw[thick,black,decorate,decoration={brace,amplitude=10pt,mirror, 
raise=5pt},xshift=0.4pt,yshift=-1.4pt](P2) -- (P5) node[black,midway,yshift=-23] {\footnotesize 
$m-1$};
\end{tikzpicture}
\]
where $m$ is odd $\ge 3$
\subsubsection*{Case~\typec{k3A}, see \cite[Theorem~1.3]{MP:cb3}}
$D\supset C$ and $\Delta(D,C)$ is as follows:
\[
\begin{tikzpicture}
 \def\sizec{0.2em}
\def\sl{0.7}
\coordinate (Q2) at (5,1*\sl);
\coordinate (P1) at (0,0);
\coordinate (P2) at (1,0);
\coordinate (P3) at (2,0);
\coordinate (P4) at (3,0);
\coordinate (P5) at (4,0);
\coordinate (P6) at (5,0);
\coordinate (P7) at (6,0);
\coordinate (P8) at (7,0);
\coordinate (P9) at (8,0);
\coordinate (P10) at (9,0);
\coordinate (P11) at (10,0);
\coordinate (Q12) at (11,0);

\draw (P2) -- (P3);
\draw [loosely dotted] (P3) -- (P4);
\draw (P4) -- (P5);
\draw (P5) -- (P6);
\draw (P6) -- (P7);
\draw (Q2) -- (P6);

\path[fill=white,draw=black] (P2) circle (\sizec)node [below, yshift=0, xshift=0] {};
\path[fill=white,draw=black] (P5) circle (\sizec)node [below, yshift=0, xshift=0] {};
\path[fill=black,draw=black] (P6) circle (\sizec)node [below, yshift=0, xshift=0] {};

\path[fill=white,draw=black] (P7) circle (\sizec)node [below, yshift=0, xshift=0] {};

\path[fill=white,draw=black] (Q2) circle (\sizec)node [below, yshift=0, xshift=0] {};

\draw[thick,black,decorate,decoration={brace,amplitude=10pt,mirror, raise=5pt},xshift=0.4pt,yshift=-1.4pt](P2) -- (P5) node[black,midway,yshift=-23] {\footnotesize $m-1$};

\end{tikzpicture}
\text{$m$ is odd $\ge 3$ \ or}
\]

\[
\begin{tikzpicture}[square/.style={regular polygon,regular polygon sides=4}]
\def\sizec{0.2em}
\def\sl{0.7}
\coordinate (1) at (0,0);
\coordinate (2) at (1,0);
\coordinate (3) at (2,0);
\coordinate (4) at (3,0);
\coordinate (5) at (4,0);
\coordinate (6) at (1,1*\sl);
\coordinate (7) at (2,1*\sl);
\coordinate (8) at (3,1*\sl);
\coordinate (9) at (2,2*\sl);
	\path (1) edge (2);
	\path (2) edge (3);
	\path (3) edge (4);
	\path (4) edge (5);
	\path (7) edge (3);
	
	\path[fill=white,draw=black] (1) circle (\sizec) ;
\path[fill=white,draw=black] (2) circle (\sizec) ;
\path[fill=black,draw=black] (3) circle (\sizec);
\path[fill=white,draw=black] (4) circle (\sizec) ;
\path[fill=white,draw=black] (5) circle (\sizec) ;
\path[fill=white, draw=black] (7) circle (\sizec);
\end{tikzpicture}
\qquad
\text{$m= 3$}
\]

\providecommand*{\BibDash}{}
\newcommand{\etalchar}[1]{$^{#1}$}
\def\cprime{$'$}

\end{document}